\documentclass[review,hidelinks,onefignum,onetabnum]{siamart220329}

\usepackage{lipsum}
\usepackage{amsfonts}
\usepackage{graphicx}
\usepackage{epstopdf}
\usepackage{algorithmic}
\ifpdf
  \DeclareGraphicsExtensions{.eps,.pdf,.png,.jpg}
\else
  \DeclareGraphicsExtensions{.eps}
\fi

\newsiamremark{remark}{Remark}
\newsiamremark{hypothesis}{Hypothesis}
\crefname{hypothesis}{Hypothesis}{Hypotheses}
\newsiamthm{claim}{Claim}

\headers{}{Y. Nabou and I. Necoara}

\title{Regularized higher-order Taylor approximation methods for composite nonlinear least-squares\thanks{Submitted to the editors \date.
\funding{ITN-ETN project TraDE-OPT funded by the European Union’s Horizon 2020 Research and Innovation Programme under the Marie Skolodowska-Curie grant agreement No. 861137;  UEFISCDI, Romania, PN-III-P4-PCE-2021-0720, under project L2O-MOC, nr. 70/2022.}}}

\author{Yassine Nabou\thanks{Automatic Control and Systems Engineering Department, National University of Science and Technology Politehnica Bucharest, Romania 
  (\email{yassine.nabou@stud.acs.upb.ro})}
\and Ion Necoara\thanks{ Automatic Control and Systems Engineering Department, National University of Science and Technology Politehnica Bucharest and Gheorghe Mihoc-Caius Iacob Institute of Mathematical Statistics and Applied Mathematics of the Romanian Academy, Bucharest, Romania (\email{ion.necoara@upb.ro})}}

\usepackage{amsopn}

\ifpdf
\hypersetup{
  pdftitle={Regularized higher-order Taylor approximation methods for nonlinear least-squares},
  pdfauthor={Y. Nabou and I. Necoara}
}
\fi

\usepackage{graphicx}
\usepackage{caption}
\usepackage{subcaption}
\providecommand{\norm}[1]{\lVert#1\rVert}
\usepackage{color}
\usepackage{amssymb}
\usepackage{amsmath}
\DeclareMathOperator*{\argmin}{arg\,min}
\newtheorem{assumption}{Assumption}
\newcommand{\red}{\textcolor{black}}
\usepackage{comment}
\usepackage{hyperref}
\usepackage{cleveref}
\usepackage{multirow}
\usepackage{adjustbox}
\nolinenumbers

\begin{document}

\maketitle

\begin{abstract}
In this paper, we develop a regularized higher-order Taylor based method for solving composite (e.g., nonlinear least-squares) problems. At each iteration, we replace each smooth component of the objective function by a higher-order Taylor approximation with an appropriate regularization, leading to a regularized higher-order Taylor approximation (RHOTA) algorithm. We derive global convergence guarantees for  RHOTA algorithm. In particular, we prove stationary point convergence guarantees for the iterates generated by RHOTA, and leveraging a Kurdyka-Lojasiewicz (KL) type property of the objective function, we derive improved rates depending on the KL parameter. When the Taylor approximation is of order $2$, we present an efficient implementation of  RHOTA algorithm, demonstrating that the resulting nonconvex subproblem can be effectively solved utilizing standard convex programming tools.  Furthermore, we extend the scope of our investigation to include the behavior and efficacy of RHOTA algorithm in handling systems of nonlinear equations and optimization problems with nonlinear equality constraints deriving new rates under improved constraint qualifications conditions. Finally, we consider solving the phase retrieval problem with a higher-order proximal point algorithm, showcasing its rapid convergence rate for this particular application. Numerical simulations on phase retrieval  and output feedback control problems also demonstrate the efficacy and performance of the proposed methods when compared to some state-of-the-art optimization methods and software.
\end{abstract}




\section{Introduction}
In this paper, we consider the following composite  problem (which, in particular, covers nonlinear least-squares):

\vspace{-0.4cm}

\begin{align}\label{ch:eq:problem}
 \min_{x \in \mathbb{R}^n} f(x): = g(F(x)) + h(x),
\end{align}

\vspace{-0.2cm}

\noindent where $F$ represents a real-vector function, defined as $F = (F_1, \ldots, F_m)$. We assume that each function $F_i:\mathbb{R}^n \to \mathbb{R}$ is $p \geq 1$ times differentiable and has the $p$th derivative Lipschitz continuous, the function $g: \mathbb{R}^m \to \mathbb{R}$ is nonsmooth, convex and Lipschitz continuous, (e.g., $g$ is the  $2$-norm) and the function $h: \mathbb{R}^n \to \Bar{\mathbb{R}}$ is proper, lower semicontinuous (lsc) and convex. Hence, $\text{dom}\, f = \text{dom}\, h$.  This formulation covers many problems from the nonlinear programming literature and appears in many real-world applications such as control, statistical estimation, grey-box minimization, machine learning, and phase retrieval \cite{DruPaq:19,DiGuMiDi:11,CaLiSo:15,NoWr:06,DuRu:19}. See also \cite{LeWr:16} for a comprehensive review on composite minimization. For example, problem  \eqref{ch:eq:problem} covers the following constrained nonlinear system of equations:
\begin{align}\label{pr:noeq}
\text{Find}\; x\in C\;\;\text{such that}\; F_i(x) = 0,\;\; i=1:m. 
\end{align}
Indeed, this problem  can equivalently be expressed as nonlinear least-squares  in the form  \eqref{ch:eq:problem}, with $h(\cdot) = 1_{C}(\cdot)$ (the indicator function of the convex set $C \subseteq \mathbb{R}^n$) and $g(\cdot) = \|\cdot\|$ (the  $2$-norm). This problem frequently arises in control \cite{DiGuMiDi:11}, phase retrieval problems  \cite{CaLiSo:15,DuRu:19,WaGiEl:17},  economic equilibrium problems \cite{DiFe:95} and learning constrained neural networks \cite{ChZu:14}.  For example, let us consider the  static output feedback stabilizability for a continuous time linear system $\dot{x} = Ax +Bu,  \; y = C x,$ where $x\in\mathbb{R}^{n_x}$ is the state vector, $u\in\mathbb{R}^{n_u\times m_u}$ is the control input and $y$ is the measured output. Using an output feedback control law of the form $u = K y$, then the system is static output feedback stabilizable if the closed loop system
$\dot{x}=A x+BK C x=(A+BK C) x$ is asymptotically stable at the origin. If the system is static output feedback stabilizable, then there always exist $X \succ 0$ and $K$ such that  $(A+BK C)^{T} X+X(A+BK C) \prec 0$ \cite{BoGhFeBa:94}. We can reformulate this bilinear matrix inequality as an equality, by introducing a matrix $Q \succ 0$  such that $(A+BK C)^{T} X+X(A+BK C)+Q=0.$ We can solve this bilinear matrix equality  by minimizing the norm of $F(X, K, Q):=(A+B K C)^{T} X+X(A+B K C) + Q$, which is a second order polynomial in $X, K$ and $Q$. Thus, the minimization problem to be solved becomes:

\vspace{-0.5cm}

\begin{align} \label{fw: eq3}
\min _{X, Q, K} \|F(X, K, Q)\|_F +  h(X,Q),
\end{align}

\vspace{-0.2cm}

\noindent where $\|\cdot\|_F$ denotes the Frobenius norm of a matrix and $h(X,Q)= \textbf{1}_{\mathbb{S}_{+}^n}(X) + \textbf{1}_{\mathbb{S}_{+}^n}(Q)$, with $\mathbb{S}_{+}^{n}$ the cone of positive definite matrices. Note that many other control problems can be posed as \eqref{fw: eq3}, see \cite{BoGhFeBa:94,Lei:16}. Another particularly interesting case of problem \eqref{ch:eq:problem} is the following optimization problem with nonlinear equality constraints:

\vspace{-0.5cm}

\begin{align}\label{eq:op_noeq}
\min_{x \in  \mathbb{R}^n} \; h(x)\;\;\;    \text{s.t.:} \;\;\;  F(x) = 0.
\end{align}

\vspace{-0.2cm}

\noindent This problem can be equivalently written using the exact penalty formulation with a given parameter $\rho > 0$ (sufficiently large) as \cite{NoWr:06}:

\vspace{-0.5cm}

\begin{align*}
    \min_{x \in \mathbb{R}^n} h(x) + \rho\|F(x)\|,
\end{align*}

\vspace{-0.2cm}

\noindent which represents a specific instance of problem \eqref{ch:eq:problem}. Optimization problems with nonlinear equality constraints appears in diverse domains including control, machine learning, signal processing and  statisticss, see, e.g., \cite{FaBo:21,Fes:20}. Specifically, the linear quadratic regulator (LQR) problem \cite{FaBo:21} is a particular case of problem \eqref{eq:op_noeq}:
\begin{align*}
    &\min_{X, K} \text{Trace}(X\Sigma) \!+\! \textbf{1}_{\mathbb{S}_{+}^n}\!(X) \;\;\; \text{s.t.}\!: \; (A \!+\! BKC)^{T}X \!+\! X(A \!+\! BKC) \!+\! C^{T} K^{T} \hat{R} K C \!+\! \hat{Q} \!=\! 0.
\end{align*}

\vspace{-0.2cm}

\noindent \textbf{Gauss-Newton type methods:} A natural approach for solving problem \eqref{ch:eq:problem} consists in linearizing the smooth part, $F$, and adding an appropriate quadratic regularization. More precisely, to obtain the next iteration, one solves the following subproblem for a given $\Bar{x}$:

\vspace{-0.5cm}

\begin{align*}
    x^{+} = \arg\min_{x \in \mathbb{R}^n} g\Big(F(\Bar{x}) + \nabla F(\Bar{x})(x - \Bar{x})\Big) +\frac{L}{2}\|x - \Bar{x}\|^2 + h(x).
\end{align*}

\vspace{-0.2cm}

\noindent This scheme, known as (proximal) Gauss-Newton method,  has been well-studied
in the literature \cite{Nes:07,DruPaq:19,NaToNe:23}. It is known that, under the Lipschitz continuity of the Jacobian, $\nabla F$, this iterative process makes the minimum norm of the subgradients of $f$ to converge to $0$ at a sublinear rate, while convergence rates under the Kurdyka-Lojasiewicz (KL) property were recently derived in \cite{Pau:16,NaToNe:23}. Trust region based Gauss-Newton methods have been also considered in \cite{CaGoTo:11} for solving problems of the form \eqref{ch:eq:problem}. The authors in \cite{CaGoTo:11} show that their proposed algorithms take at most $\mathcal{O}(\epsilon^{-2})$ function evaluations to reduce the size of a first-order criticality measure below a given accuracy $\epsilon$. While these schemes have shown empirical success in addressing challenging and complex optimization problems, their convergence rates are known to be slow. In order to speed up the convergence rates, one needs to use higher-order information (derivatives) to construct a more accurate higher-order (Taylor) model that effectively approximates the objective function. 
\noindent From our knowledge, there are very few studies considering the utilization of higher-order derivatives to address problems of the form \eqref{ch:eq:problem} where the function $g$ is both convex and Lipschitz continuous. For example, in \cite{CaGoTo:13}, the authors explore the scenario where $g(\cdot) = \frac{1}{2}\|\cdot\|^2$, $h(\cdot) = 0$ and the next iterate, given the current pointy $\Bar{x}$, is the minimizer of the following cubic subproblem:

\vspace{-0.5cm}

\begin{align*}
x^{+} = \arg\min_{x \in \mathbb{R}^n }\;  (x - \Bar{x})^{T}J_F(\Bar{x})^{T}F(\Bar{x}) + \frac12(x - \Bar{x})^T \Bar{B}(\Bar{x}) (x - \Bar{x}) + \frac{M}{3} \|x - \Bar{x}\|^3,
\end{align*}

\vspace{-0.1cm}

\noindent where $M >0$, $J_F(\Bar{x})$ is the Jacobian of $F$ at $\Bar{x}$ and $\Bar{B}(\Bar{x})$ is an approximation of the true hessian of the function $\frac{1}{2}\|F(x)\|^2$ at $\Bar{x}$. When the residuals $F_i$, the Jacobian $\nabla F$ and the Hessian $\nabla^2 F_i$ for each $i\in\{1,\cdots,m\}$ are simultaneuosly  Lipschitz continuous on a neighborhood of $\bar{x}$,  \cite{CaGoTo:13} shows that this scheme takes at most $\mathcal{O}\left(\epsilon^{-\frac{3}{2}}\right)$ residual and Jacobian evaluations to drive either the Euclidean norm of the residual or its gradient below $\epsilon$. Further, in \cite{GoReSc:19} a similar approach is adopted, with $g(\cdot) = \frac{1}{2}\|\cdot\|^2$ and  $h(\cdot) = 0$, by constructing a  quadratic approximation of $F$  alongside an appropriate regularization, i.e., given the current point $\Bar{x}$, the next iterate is the minimizer of the following $r$-regularized subproblem:

\vspace{-0.5cm}

\begin{align*}
\min_{x\in \mathbb{R}^n }\frac{1}{2} \left\|F(\Bar{x}) + \nabla F(\Bar{x})(x - \Bar{x}) + \frac{1}{2} (x - \Bar{x})^T\nabla^2 F(\Bar{x})(x - \Bar{x})\right\|^2\!
     + \frac{M}{r}\|x - \Bar{x}\|^r,
\end{align*}

\vspace{-0.1cm}

\noindent where $r\geq 2$ is a given constant and $M>0$ is a regularization parameter. Paper \cite{GoReSc:19} establishes convergence to an $\epsilon$ first-order stationary point of the objective within $\mathcal{O}\left( \epsilon^{-\min\left(\frac{r}{r-1}, \frac{3}{2}\right)} \right)$ iterations, provided that the residuals $F_i$, the Jacobian $\nabla F$ and the Hessian $\nabla^2 F_i$  for each $i\in\{1,\cdots,m\}$ are  Lipschitz continuous on a neighborhood of a stationary point.  It's important to note that the aforementioned subproblem is at least quartic and thus hard to solve. Therefore, using the norm $\|\cdot\|$ as the merit function instead of $\|\cdot\|^2$ is more beneficial, since in the later case the condition number usually doubles, although the objective function for the first choice is  nondifferentiable \cite{Nes:07}. This strategy has been explored in \cite{CaGoTo:22}, wherein the authors introduce an adaptive higher-order trust-region algorithm for solving problem \eqref{ch:eq:problem} with $h$ smooth, where $F$ and $h$ are approximated  with higher-order Taylor expansions. Paper \cite{CaGoTo:22} establishes convergence of order $\mathcal{O}\left(\epsilon^{-\frac{p+1}{p}}\right)$ to achieve a reduction in a given criticality measure below a prescribed accuracy $\epsilon$. Note that, the optimization problem, the algorithm, and consequently the convergence analysis in \cite{CaGoTo:22} are different from the present paper. Moreover, it remains open  whether one can solve the corresponding  subproblem in \cite{CaGoTo:22} efficiently  for $p\geq 2$, along with establishing convergence rates under the Kurdyka-Lojasiewic (KL) property. 

\medskip

\noindent \textbf{Contributions:} In this paper, we present a regularized higher-order Taylor approximation method (called RHOTA) for solving problem \eqref{ch:eq:problem}. At each iteration of RHOTA, we build a higher-order composite model and minimize it to compute the next iteration. We also derive  worst-case complexity bounds for the RHOTA algorithm. Thus, our main contributions can be summarized as follows:

\vspace{0.1cm}

\noindent (i) We present a \textit{new regularized higher-order Taylor approximation  algorithm} (called RHOTA) for solving  \eqref{ch:eq:problem}. At each iteration, we replace the smooth components of $F$ with  higher-order Taylor approximations of order $p \geq 1$ and add a proper regularization. The  minimizer of this approximate model yields the next iteration. An adaptive  variant of RHOTA is also presented.

\vspace{0.1cm}

\noindent (ii) We derive convergence guarantees for RHOTA algorithm under different assumptions on the problem \eqref{ch:eq:problem} and  \textit{two optimality measures} characterizing first-order stationary points. More precisely, we show that the iterates generated by RHOTA converge globally to near-stationary points and the convergence rate is of order $\mathcal{O}\left(\epsilon^{-\frac{p+1}{p}}\right)$, where $\epsilon$ is the desired accuracy. Additionally, when the objective function satisfies a Kurdyka-Łojasiewicz (KL) type property, we derive linear/sublinear convergence rates in function values, depending on the parameter of the KL condition.

\vspace{0.1cm}

\noindent (iii) When $g(\cdot) = \|\cdot\|$ and $h(\cdot)$ is quadratic, we present an efficient implementation of RHOTA algorithm for $p=2$. In particular, we show that the resulting nonconvex subproblem is equivalent to minimizing an explicitly written convex function over a convex set that can be solved using standard convex tools. This represents a significant advancement towards the practical implementation of higher-order methods for solving composite (e.g., nonlinear least-squares) problem \eqref{ch:eq:problem}.

\vspace{0.1cm}

\noindent (iv) We also analyze  the convergence  behavior of RHOTA algorithm when employed to address systems of nonlinear equations and optimization problems featuring nonlinear equality constraints and derive  new global convergence rates for our algorithm on these classes of problems under \textit{improved constraint qualification conditions}. 

\vspace{0.1cm}

\noindent (v) Finally, we demonstrate that two important applications, the phase retrieval \cite{CaLiSo:15,DuRu:19,WaGiEl:17} and the output feedback control \cite{DiGuMiDi:11,Lei:16}, can be effectively framed within the context of  problem  \eqref{ch:eq:problem}. For phase retrieval, RHOTA algorithm yields  a higher-order proximal point algorithm (called HOPP) and  our analysis establishes rapid convergence of HOPP on this application. RHOTA also yields an efficient regularized Gauss-Newton  algorithm for solving the  output feedback control problem.  Numerical results on output feedback control for linear systems from  COMPl$_{e}$ib library \cite{Lei:16}  and on image recovery on handwritten digit images from MINIST library  \cite{MNIS}  demonstrate  the superior performance of  RHOTA  compared to some state-of-the-art optimization method \cite{DuRu:19} and software \cite{DiGuMiDi:11} developed specifically for these applications. 

\vspace{0.1cm}

\noindent \textit{Content}.  The remainder of this paper is organized as follows: in Section 2, we present our settings; in Section 3, we introduce our main assumptions and  RHOTA algorithm; In Section 4 we  derive the main convergence results;  in Section 5, we present an efficient implementation of RHOTA  for $p=2$; in Section 6, we focus on the convergence behaviour of RHOTA when solving  specific problems from nonlinear programming; finally, in Section 7, we illustrate the effectiveness of the proposed algorithms within the context of phase retrieval and output feedback control problems using real data.



\vspace{-0.1cm}

\section{Notations and preliminaries}
We consider a  finite-dimensional Euclidean space $\mathbb{R}^n$  endowed with an  inner product $\langle s,x\rangle$ and the corresponding norm $\|s\| =\langle s,s\rangle^{1/2}$ for any $s,x\in  \mathbb{R}^n$. For a twice differentiable function $\phi$ on a convex and open domain $\text{dom}\,\phi \subseteq \mathbb{R}^n$, we denote by $\nabla \phi(x)$ and $\nabla^{2} \phi(x)$ its gradient and hessian evaluated at $x\in \text{dom}\,\phi$, respectively.  Throughout the paper, we consider $p$ a positive integer.  In what follows, we often work with directional derivatives of function $\phi$ at $x$ along directions $h_{i}\in \mathbb{R}^n$ of order $p$, $D^{p} \phi (x)[h_{1},\cdots,h_{p}]$, with $i=1:p$.  If all the directions $h_{1},\cdots,h_{p}$ are the same, we use the notation
$ D^{p} \phi(x)[h]$, for $h\in\mathbb{R}^n$.  Note that if $\phi$ is $p$ times differentiable, then $D^{p} \phi (x)$ is a symmetric $p$-linear form and its norm is defined as \cite{Nes:20}: $\norm{D^{p} \phi (x)}= \max\limits_{h\in\mathbb{R}^n} \left\lbrace D^{p} \phi (x)[h]^{p} :\norm{h}\leq 1 \right\rbrace.$
\noindent Further, the Taylor approximation of order $p$ of $\phi$ at $x\in \text{dom}\,\phi$ is denoted with:

\vspace{-0.6cm}

\begin{align*}
T_p^{\phi}(y; x)= \phi(x) + \sum_{i=1}^{p} \frac{1}{i!} D^{i} \phi(x)[y-x]^{i} \quad \forall y \in \mathbb{R}^n.
\end{align*}

\vspace{-0.3cm}

\noindent Let $\phi: \mathbb{R}^n \mapsto\bar{\mathbb{R}}$ be a $p$ differentiable function on $\text{dom}\,\phi$. Then, the $p$ derivative is Lipschitz continuous if there exist a constant $L_p^{\phi} > 0$ such that:

\vspace{-0.5cm}

   	\begin{equation} \label{eq:001}
    	\| D^p \phi(x) - D^p \phi(y) \| \leq L_p^{\phi} \| x - y \| \quad  \forall x,y \in \text{dom}\,\phi.
    	\end{equation}  

\vspace{-0.2cm}
     
\noindent It is well known that if \eqref{eq:001} holds,  then the residual between the function  and its Taylor approximation can be bounded \cite{Nes:20}:

\vspace{-0.4cm}

    \begin{equation}\label{eq:TayAppBound}
    \vert \phi(y) - T_p^{\phi}(y;x) \vert \leq  \frac{L_p^{\phi}}{(p+1)!} \norm{y-x}^{p+1} \quad  \forall x,y\in \text{dom}\,\phi.
    \end{equation}



\vspace{-0.3cm}

\noindent Let set $\Omega\subseteq \mathbb{R}^n$. The Fr$\acute{e}$chet regular normal cone to $\Omega$ at $\bar{x} \in\Omega$ is defined by \cite{Roc:70}: 

\vspace{-0.5cm}

\begin{align*}
    \widehat{\mathcal{N}}_{\Omega}(\Bar{x}) = \Big \{w\in\mathbb{R}^n:\; \limsup_{x \xrightarrow[\Omega]{} \bar{x}} \frac{\langle w,x - \bar{x} \rangle}{\|x - \bar{x}\|} \leq 0 \Big \},
\end{align*}

\vspace{-0.3cm}

\noindent 
The limiting normal cone to $\Omega$ at $\bar{x}$ is defined as~\cite{Roc:70}:

\vspace{-0.3cm}

\begin{align*}
    \mathcal{N}_{\Omega}(\bar{x}) = \left\{w\in\mathbb{R}^n:\; \exists \bar{x}_k \xrightarrow[\Omega]{} \bar{x}, w_k\to w\; \text{as}\; k\to \infty\; \text{with}\; w_k\in \widehat{\mathcal{N}}_{\Omega}(x_k)\right\}.
\end{align*}

\vspace{-0.2cm}

\noindent The epigraph of a proper lower semicontinuous (lsc) function $\phi: \mathbb{R}^n \to \bar{\mathbb{R}}$ is $\text{epi}\,\phi = \{(x,\alpha)\in\mathbb{R}^n\times\mathbb{R}: \phi(x)\leq \alpha\}$. A function $\phi:\mathbb{R}^n \to \bar{\mathbb{R}}$ is called regular at $\bar{x}$ if $\phi(\bar{x})$ is finite and $\text{epi}\,\phi$ is Clarke regular at $(\bar{x},\phi(\bar{x}))$, i.e., $\text{epi}\,\phi $ is locally closed and it holds that $\mathcal{N}_{\text{epi}\,\phi}(\bar{x},\phi(\bar{x})) = \widehat{\mathcal{N}}_{\text{epi}\,\phi}(\bar{x},\phi(\bar{x}))$ (see Definition 7.25 in \cite{Roc:70}). For any  $x\in\text{dom}\, \phi$,  the limiting subdifferential and the horizon subdifferential of a proper lsc function $\phi$ at $x$ can be characterized via the limiting normal cone (see Theorem 8.9 in \cite{Roc:70}):

\vspace{-0.4cm}

\begin{align*}
\partial \phi (x) = \left\{g_x \! : (g_x ,-1) \in \mathcal{N}_{\text{epi}\,\phi}(x,\phi(x))    \right\}, \;\;\;    \partial^{\infty} \phi (x) = \left\{g_x \! : (g_x ,0) \in \mathcal{N}_{\text{epi}\,\phi}(x,\phi(x))\right\}. 
\end{align*}

\vspace{-0.1cm}

\noindent Obviously, $\partial^{\infty} \phi (x)$ is a cone. Denote  $S_\phi(x) := \text{dist}(0,\partial \phi(x))$ the  minimum norm subdifferential. For a given $x_0\in\text{dom}\,\phi$, we denote by $ \mathcal{L}_\phi(x_0)$ the level set of the function $\phi$, i.e., $ \mathcal{L}_\phi(x_0) := \{x\in\text{dom}\,\phi: \; \phi(x) \leq \phi(x_0)\}$. Moreover, if $g$ is proper lsc convex and Lipschitz continuous function and $F = (F_1,\cdots,F_m)$, where $F_i$'s are differentiable functions, then the following chain rule applies at any $x \in  \text{dom}\,g(F)$ (see Theorem 10.6 in \cite{Roc:70}): $\partial (g\circ F)(x) \subseteq \nabla F(x)^{T}\partial g(F(x))$. The equality holds if, in addition, $g$ is regular at $F(x)$. Further, given $f = f_1 + f_2$, where $f_i$'s are proper lsc  functions either both convex or $f_1$ locally Lipschitz, then for all $x\in\text{dom}\, f$ we have (see Corollary 10.9 in \cite{Roc:70}):  $\partial f(x)\subseteq \partial f_1(x) + \partial f_2(x)$. The equality holds e.g., if each convex $f_i$, for $i=1:2$, is also regular at $x$ or if $f_1$ is strictly differentiable. Hence, according to these calculus rules, for the objective function in \eqref{ch:eq:problem}  if $g$ is proper lsc convex function and regular at $F(x)$ and $h$ is locally Lipschitz, then we have:

\vspace{-0.5cm}

\begin{align}
\label{chain-rule}
\partial f(x) \subseteq \nabla F(x)^T \partial g(F(x)) + \partial h(x),
\end{align}

\vspace{-0.2cm}

\noindent while equality holds under additional conditions on $h$. For a given function $f$, any point $x^*$ such that $0\in\partial f(x^*)$ is called \textit{stationary point} and denote $\texttt{stat} f := \{x^*: 0 \in \partial f(x^*) \}$. When dealing with functions of the form \eqref{ch:eq:problem}, different types of stationarity have been considered in the literature. Clearly, the first choice is based on  points $x^*$ such that $0\in\partial f(x^*)$. A second choice, considered e.g., in \cite{Yua:85,CaGoTo:13}, is based on points $x^*$ for which the  directional derivative satisfies: 

\vspace{-0.5cm}

\begin{align}\label{eq:stat_yan}
D f(x^*;d):=\sup_{\lambda\in\partial g(F(x^*)),\bar{\lambda}\in\partial h(x^*)}  \left(\nabla F(x^*)^T \lambda + \bar{\lambda}\right)^T d  \geq 0 \;\;\;\;   \forall d\in\mathbb{R}^n.    
\end{align}

\vspace{-0.2cm}

\noindent The stationarity condition \eqref{eq:stat_yan} is  equivalent to the first one (i.e.,  $0\in\partial f(x^*)$) if e.g., \eqref{chain-rule} holds with equality and $g, h$ are  locally Lipschitz. Indeed, if \eqref{eq:stat_yan} holds and $g,h$ are locally Lipschitz, then there exist finite $\lambda\in\partial g(F(x^*)) $ and $\bar{\lambda} \in\partial h(x^*)$ such that:
\begin{align*}
\left(\nabla F(x^*)^T \lambda + \bar{\lambda}\right)^Td \geq 0 \;\;\;\; \forall d\in\mathbb{R}^n. 
\end{align*}
Choosing $d = - (\nabla F(x^*)^T \lambda + \bar{\lambda})$, we get that $0 = \nabla F(x^*)^T \lambda + \bar{\lambda}\in \partial f(x^*)$, provided that  \eqref{chain-rule} holds with equality. For the other implication, let us  assume that $0\in\partial f(x^*)$.  Then, if \eqref{chain-rule} holds, it follows that there exist finite $\lambda\in\partial g(F(x^*)) $ and $\bar{\lambda} \in\partial h(x^*)$ such that  $\nabla F(x^*)^T \lambda + \bar{\lambda} = 0$. Hence, we get:

\vspace{-0.5cm}

\begin{align*}
 D f(x^*;d) \geq  \left(\nabla F(x^*)^T \lambda + \bar{\lambda}\right)^T d = 0 \; \;\;\forall d\in\mathbb{R}^n, \;\; \text{i.e., eq.} \;  \eqref{eq:stat_yan}.
\end{align*}

\vspace{-0.2cm}

\noindent Any point $x^*$ satisfying \eqref{eq:stat_yan} is called a \textit{weak stationary point} (note that if $x^*$ satisfies $0\in\partial f(x^*)$, then from \eqref{chain-rule} there exist $\lambda\in\partial g(F(x^*))$ and $\bar{\lambda}\in\partial h(x^*)$ such that $\nabla F(x^*)^T \lambda + \bar{\lambda} = 0$ and thus \eqref{eq:stat_yan} also holds at $x^*$; on the other hand, if $x^*$ satisfies \eqref{eq:stat_yan}, then one needs additional assumptions in order to also have $0\in\partial f(x^*)$, such as boundedness of $\partial h(x^*)$ and  the chain rule \eqref{chain-rule} must hold with equality). Clearly, any local minimum $x^*$ of problem \eqref{ch:eq:problem} satisfies the two previous stationary point conditions.
\noindent Next, we recall the \textit{Kurdyka-Lojasiewicz (KL)} property for semi-algebraic functions  around a compact set~$\Omega$ \cite{BolSab:14}:

\vspace{-0.5cm}

\begin{equation}\label{eq:kl}
\phi(x) - \phi_*  \leq \sigma_q  S_{\phi}(x)^q \quad   \forall x\!: \;  \text{dist}(x, \Omega) \leq \delta, \; \phi_* < \phi(x) < \phi_* + \epsilon.  
\end{equation}

\vspace{-0.2cm}




\noindent  The relevant aspect of the KL property is when $\Omega$ is a subset of stationary points for $\phi$, i.e.  $\Omega \subseteq \texttt{stat} \phi$, since it is easy to establish the KL property when $\Omega$ is not related to stationary points. Note that the set of semi-algebraic functions include real polynomial functions, vector or matrix (semi)norms (e.g., $\|\cdot\|_p$ with $p \geq 0$ rational number),  uniformly convex functions, as well as composition of semi-algebraic functions (see \cite{BolSab:14} for a comprehensive list).


\section{Regularized higher-order Taylor approximation method}
\label{sec:3}
In this section, we present a regularized higher-order Taylor approximation algorithm for solving composite problem \eqref{ch:eq:problem}. We consider the following assumptions:
\begin{assumption}
\label{as:1}
The following statements hold for  optimization problem \eqref{ch:eq:problem}:
\begin{enumerate}
    \item For $F = (F_1,\cdots,F_m)$, each component  $F_i$ is $p$ times differentiable function  with the $p$th derivative Lipschitz continuous with constant $L_p^i$.
    \item Function $g$ is convex, Lipschitz continuous with constant $L_g$ and $h$ is proper lower semicontinuous and simple convex function.
    \item Problem \eqref{ch:eq:problem} has a solution and hence $\inf_{x \in \emph{dom}\;f}f(x) \geq f^*$.
\end{enumerate}
\end{assumption}

\noindent From Assumption \ref{as:1} and the inequality \eqref{eq:TayAppBound}, we get for all $i=1:m$: 

\vspace{-0.5cm}

\begin{align}
    \left |F_i(x) - T_p^{F_i}(x;y) \right| \leq \frac{L_p^i}{(p+1)!}\|y - x\|^{p+1}\;\;\;  \forall x,y\in\mathbb{R}^n. 
\end{align}

\vspace{-0.2cm}

\noindent Further, using that the function $g$ is Lipschitz continuous, we get the following inequality valid for all $x,y\in\mathbb{R}^n$:

\vspace{-0.5cm}

\begin{align}\label{main:eq}
    \left|g(F(x)) - g\left(T_p^F(x;y)\right) \right|\leq L_g\left\| F(x) -  T_p^F(x;y)\right\| \leq\frac{L_g\|L_p\|}{(p+1)!}\|x - y\|^{p+1},
\end{align}

\vspace{-0.2cm}

\noindent where $T_p^F(x;y) = \Big(T_p^{F_1}(x;y),\cdots,T_p^{F_m}(x;y)\Big)$ and $L_p = \left(L_p^1,\cdots,L_p^m\right)$. Then, based on this upper bound approximation of the objective function, one can consider an iterative process, where given the current iterate, $\Bar{x}$, and  a proper regularization parameter $M>0$, the next point is computed from the following subproblem:

\vspace{-0.5cm}

\begin{align}\label{sub:1}
 x \leftarrow \arg\min_{y \in \mathbb{R}^n} s_M(y;\bar{x}):= g\left(T_p^F(y;\bar{x})\right) + \frac{M}{(p+1)!}\|y - \bar{x}\|^{p+1} + h(y).      
\end{align}

\vspace{-0.2cm}

\noindent Note that if $x = \Bar{x}$ in the previous subproblem, then $x$ is a stationary point of the original problem \eqref{ch:eq:problem}. Note also that for $p=1$, this algorithm reduces to the regularized Gauss-Newton method analyzed in \cite{Nes:07,DruPaq:19,NaToNe:23}. Now we are ready to present our regularized higher-order Taylor approximation method, called RHOTA (see Algorithm \ref{alg:RHOTA}).
\begin{algorithm}
\begin{algorithmic}
\STATE{Given $x_0 \in \text{dom}\, f$ and $ M>0$. \\
For $k \geq 0$ do:}
\STATE{compute $x_{k+1} \in \text{dom}\,f$ \textit{inexact} solution of subproblem \eqref{sub:1} satisfying the  \textit{descent}:}
\begin{align}\label{des:1}
  s_M(x_{k+1};x_k)\leq s_M(x_k;x_k).
\end{align}
\end{algorithmic}
\caption{RHOTA}
\label{alg:RHOTA}
\vspace{-0.2cm}
\end{algorithm}
\noindent Note that usually the subproblem \eqref{sub:1} is nonconvex for any $p\geq 2$. In order to get  descent for the sequence $(f(x_k))_{k\geq 0}$, it is enough to assume that $x_{k+1}$ satisfies the descent \eqref{des:1}. However, to derive convergence rates to stationary points or in function values (under the KL property), we need to require additionally properties for $x_{k+1}$, e.g., $x_{k+1}$ generated by algorithm must satisfy some inexact (local) optimality condition as in \eqref{eq:inx} (see Theorem \ref{th:2}) or as in \eqref{eq:inx_cr} (see Theorem \ref{th4.6}), i.e., computing a  minimizer of the Taylor based model $s_M(\cdot;x_k)$ within an Euclidean ball. We show in Section \ref{sec:sol} that one can still use the powerful tools from convex optimization to solve the \textit{nonconvex} subproblem \eqref{sub:1} globally for some particular choices of $p>1$. More precisely, when the outer function $g$ is the norm and the Taylor approximation is of order $p=2$, we show that the corresponding subproblem can be solved globally by efficient convex  algorithms.


\section{Convergence analysis of RHOTA}
In this section, we analyze the convergence behavior of RHOTA algorithm under different assumptions for problem \eqref{ch:eq:problem}, i.e., under Assumption \ref{as:1}  and,  additionally, the objective function satisfies the KL condition. First, let us prove that the sequence $(f(x_k))_{k\geq 0}$ is a nonincreasing sequence.

\begin{theorem}\label{th:1}
Let Assumption \ref{as:1} hold and let $(x_k)_{k\geq 0}$ be generated by RHOTA with $M - L_g\|L_{p}\| > 0$. Then, we have:

\noindent 1. The sequence $(f(x_k))_{k\geq 0}$ is nonincreasing and satisfies:

\vspace{-0.3cm}

\begin{align}\label{eq:dec}
f(x_{k+1})  \leq f(x_k) - \frac{M - L_g\|L_{p}\|}{(p+1)!} \|x_{k+1} - x_k \|^{p+1}. 
\end{align}

\vspace{-0.1cm}

\noindent 2. The sequence $(x_k)_{k\geq 0}$ satisfies:

\vspace{-0.5cm}

\begin{align*}
&\sum_{i=1}^{\infty} \|x_{k+1} \!- x_k\|^{p+1} \!< \infty,\;\; \lim_{k\to \infty} \|x_{k+1} \!-\! x_k\| =\! 0 \;\; \emph{and}\;\; \min_{j=0:k}\|x_{j+1} \!- x_j\|^{p+1} \leq \! \mathcal{O}\!\left(\frac{1}{k}\right).
\end{align*}

\vspace{-0.1cm}

\end{theorem} 
\begin{proof}
From inequality \eqref{main:eq}, we get:
\begin{align*}
    -\frac{L_g\|L_{p}\|}{(p+1)!}\|x_{k+1} - x_k\|^{p+1} + g(F(x_{k+1}))\leq g\left(T_{p}^F(x_{k+1};x_k)\right).
\end{align*}
Further, using the descent \eqref{des:1}, we also get:
\begin{align*}
\frac{M - L_g\|L_{p}\|}{(p+1)!}&\| x_{k+1} - x_k\|^{p+1} + f(x_{k+1})\\
&\leq g\left(T_{p}^F (x_{k+1};x_k)\right) + h(x_{k+1}) + \frac{M}{(p+1)!}\|x_{k+1} - x_k\|^{p+1}\\
& = s_M(x_{k+1};x_k)\leq s_M( x_k ;x_k) = f(x_k).
\end{align*}
Hence, the sequence $(f(x_k))_{k\geq 0}$ is monotonically nonincreasing. Further, summing up the last inequality and using that $f$ is bounded from below by $f^*$, we get:
\begin{align*}
\sum_{j=0}^{k}\frac{M - L_g\|L_{p}\|}{(p+1)!}\|x_{j+1} - x_j\|^{p+1} 
    \leq f(x_0) - f(x_k) \leq f(x_0) - f^*.
\end{align*}
Hence, there exists  $\bar{k} \in\{0,\cdots,k\}$ such that:
\begin{align}\label{eq:cv_pnt}
   \|x_{\bar{k}+1} - x_{\bar{k}}\|^{p+1} =  \min_{j = 0:k}\|x_{j+1} - x_j\|^{p+1} \leq \frac{(f(x_0) - f^*)(p+1)!}{(M-L_g\|L_p\|)(k+1)},
\end{align}
and then our assertions follow.
\end{proof}
\begin{remark}
Theorem \ref{th:1} requires $M - L_g\|L_{p}\| > 0$, where $\|L_{p}\|\!=\!\| (L_{p}^1,\!\cdots\!,L_{p}^m)\|$. If $L_g$ and $(L_{p})_{i=1}^m$ are known, then one can choose $M = L_g\|L_{p}\| + R_0$ for some $R_0>0$. In Section \ref{sec:ad}, we propose an adaptive variant of RHOTA that does not require the knowledge of the Lipschitz constants $L_g$ and $(L_{p})_{i=1}^m$. 
\end{remark}


\subsection{Global convergence analysis and rates}\label{sec:subg}
\noindent In this section, we derive global convergence rates for the iterates of RHOTA algorithm in two  optimality measures: minimum norm subdifferential and the optimality measure introduced in \cite{CaGoTo:13,CaGoTo:22,Yua:85}.

\subsubsection{Minimum norm subgradients convergence rate}
In this section, \textit{we derive a global first-order convergence rate for RHOTA  to stationary points}, i.e., we use the minimum norm subdifferential as a measure of optimality. In \cite{DruPaq:19,NabNec:22}, the authors prove for  composite problem \eqref{ch:eq:problem} that the quantity $\text{dist}(0,\partial f(x_{k+1}))$ doesn't invariably approach zero as $\|x_{k+1} - x_k\|$ tends to zero. Hence, in alignment with the framework outlined in \cite{NabNec:22}, for a given $\mu>0$ we introduce the (artificial) sequence:
\begin{align}\label{eq:tk1}
    y_{k+1} = \argmin_{y \in \mathbb{R}^n} f(y) + \frac{\mu}{(p+1)!}\|y - x_k\|^{p+1}.
\end{align}
From the optimality conditions of the iteration $y_{k+1}$, we get:

\vspace{-0.5cm}

\begin{align*}
-\frac{\mu}{p!}\|y_{k+1} - x_k\|^{p-1}\left(y_{k+1} - x_k\right) \in \partial f(y_{k+1}).
\end{align*}

\vspace{-0.3cm}

\noindent This implies that

\vspace{-0.5cm}

\begin{align}\label{eq:seq-y}
 \mathcal{S}_f(y_{k+1})=\text{dist}\left(0,\partial f(y_{k+1})\right) \leq \frac{\mu}{p!}\|y_{k+1} - x_k\|^p.
\end{align}

\vspace{-0.2cm}

\noindent In the next theorem we establish that the sequence $(y_{k})_{k\geq 0}$ is close to the sequence $(x_{k})_{k\geq 0}$, both sequences have the same set of limit points, and $(y_{k})_{k\geq 0}$ converges  towards a stationary point of the original problem with a rate $\mathcal{O}(k^{-\frac{p}{p+1}})$. These results are  valid under the condition that the sequence $(x_{k})_{k\geq 0}$ generated by  RHOTA algorithm satisfies an inexact optimality criterion.

\begin{theorem}\label{th:2}
Let the assumptions of Theorem \ref{th:1} hold. Let $(x_k)_{k\geq 0}$ be generated by RHOTA algorithm. Let $\mu>M + L_g\|L_{p}\|$ and $y_{k+1}$ be given in \eqref{eq:tk1} and assume  $x_{k+1}$ satisfies the following inexact optimality condition for subproblem \eqref{sub:1}:
\begin{align}\label{eq:inx}
s_M(x_{k+1} ; x_k)  - \min_{ {y : \; \|y - x_k\|\leq D_k }} s_M(y;x_k) \leq \frac{\delta}{(p+1)!}\norm{x_{k+1} - x_k}^{p+1},
\end{align}

\vspace{-0.3cm}

\noindent where $\delta \geq 0$, $D_k:= \left(\frac{(p+1)!}{\mu}(f(x_k) - f^*)\right)^{\frac{1}{p+1}}$ and $s_M(\cdot;x_k)$ is given in \eqref{sub:1}. If we denote $L_\mu  =  \left(\frac{\mu + \delta+ L_g\|L_p\| - M}{\mu - (M + L_g\|L_{p}\|)}\right)$, then  we have:
\begin{enumerate}
    \item The sequences $(y_{k})_{k\geq 0}$  satisfies   $\| y_{k+1} - x_k \| ^{p+1}\leq L_\mu \|x_{k+1} - x_k\|^{p+1}$ \quad $\forall k \geq 0$.
    \item  The following convergence rate holds:
    \begin{align*}
   \min_{j=0:k}S_f(y_{j+1})^{\frac{p+1}{p}} \leq\frac{L_{\mu}(f(x_0) -  f^*)}{k+1}\left(\frac{\mu}{p!}\right)^{\frac{p+1}{p}}\frac{(p+1)!}{M - L_g\|L_{p}\|}. 
    \end{align*}
\end{enumerate}
\end{theorem}
\begin{proof}
From the definition of $y_{k+1}$, we have:
\begin{align*}
 &f(y_{k+1}) +\frac{\mu}{(p+1)!}\|y_{k+1} - x_k\|^{p+1} \stackrel{\eqref{eq:tk1}}{\leq} f(x_{k+1}) + \frac{\mu}{(p+1)!}\|x_{k+1} - x_k \|^{p+1}\\
 & \stackrel{\eqref{main:eq}}{\leq} s_M(x_{k+1};x_k) + \frac{\mu + L_g\|L_{p}\|-M}{(p+1)!}\|x_{k+1} - x_k\|^{p+1}\\
 &\stackrel{\eqref{eq:inx}}{\leq}\!\!\! \min_{ y: \; \| y - x_k\|\leq D_k }\!\!\!\! s_M(y;x_k) \!+\! \frac{\delta}{(p+1)!} \norm{x_{k+1} \!-\! x_k}^{p+1} \!+\! \frac{\mu+ L_g\|L_{p}\|-M}{(p+1)!}\|x_{k+1} \!-\! x_k\|^{p+1}\\ 
 &\stackrel{\eqref{main:eq}}{\leq}\!\!\!\! \min_{ y: \; \|y - x_k\|\leq D_k }\!\! f(y) + \frac{M \!+\! L_g\|L_{p}\|}{(p\!+\!1)!} \|y - x_k\|^{p+1} \!+\! \frac{\mu \!+\! \delta+ L_g\|L_{p}\|\!-\!M}{(p\!+\!1)!}\|x_{k+1} \!-\! x_k\|^{p+1}\\
 &\leq f(y_{k+1}) + \frac{M + L_g\|L_{p}\|}{(p+1)!} \|y_{k+1} - x_k\|^{p+1}  + \frac{\mu +\delta + L_g\|L_p\| - M}{(p+1)!}\|x_{k+1} - x_k\|^{p+1},
\end{align*}
where the last inequality is derived from the observation that $\|y_{k+1} - x_k\|\leq D_k$. Indeed, from the definition of $y_{k+1}$ from \eqref{eq:tk1} we have:
\begin{align*}
f(y_{k+1}) + \frac{\mu}{(p+1)!}\|y_{k+1} - x_k\|^{p+1}\leq f(x_k),    
\end{align*}
which implies that:
\begin{align*}
    \frac{\mu}{(p+1)!}\| y_{k+1} - x_k \|^{p+1}\leq f(x_k) - f^*.
\end{align*}
and thus $\|y_{k+1} - x_k\|\leq D_k$. Further, we have:
\begin{align}
\label{eq:y-x}
    \|y_{k+1} - x_k\|^{p+1} &\leq \left(\frac{\mu + \delta + L_g\|L_p\| - M}{\mu - (M + L_g\|L_{p}\|)}\right)\|x_{k+1} - x_k\|^{p+1} \quad \forall k \geq 0,
\end{align} 
which is the first assertion. It follows immediately from the last inequality that $(y_k)_{k\geq 0}$ and $(x_k)_{k\geq 0}$ have the same set of limit points. Additionally, from \eqref{eq:cv_pnt}, we have that there exists  $\bar{k} \in\{0,\cdots,k\}$ such that:
\begin{align}\label{eq:bn_sub}
  & \min_{j=0:k}S_f(y_{j+1})^{\frac{p+1}{p}}   \leq  S_f(y_{\bar{k} + 1})^{\frac{p+1}{p}} \stackrel{\eqref{eq:seq-y}}{\leq} \left(\frac{\mu}{p!}\right)^{\frac{p+1}{p}}\|y_{\bar{k}+1} - x_{\bar{k}}\|^{p+1} \\
  & \stackrel{\eqref{eq:y-x}}{\leq} L_{\mu}\left(\frac{\mu}{p!}\right)^{\frac{p+1}{p}}\|x_{\bar{k}+1} - x_{\bar{k}}\|^{p+1} \stackrel{\eqref{eq:cv_pnt}}{\leq} L_{\mu}\left(\frac{\mu}{p!}\right)^{\frac{p+1}{p}}\frac{(f(x_0) - f^*)(p+1)!}{(M - L_g\|L_{p}\|)(k+1)}.\nonumber
\end{align}
Hence, our second statement follows.  
\end{proof}
\begin{remark}
In Theorem \ref{th:2} we establish convergence rate guarantees of order $\mathcal{O}\left(k^{-\frac{p}{p+1}}\right)$ to  (near) stationary points, which is the usual convergence rate for higher-order algorithms for (unconstrained) nonconvex $p$-smooth problems \cite{BiGaMa:17,CaGoTo:22,NabNec:22}. In our convergence analysis, we additionally assume that the sequence generated by  RHOTA algorithm satisfy an inexact optimality condition \eqref{eq:inx}, which requires computing a  minimum point  over an Euclidean ball. In Section \ref{sec:sol} we show that we can compute such a point in RHOTA algorithm for the particular case  $g(\cdot)=\|\cdot\|$ and $p=2$. More precisely, we show that one can still use powerful tools from convex optimization to even compute the \textit{global solution of the nonconvex} subproblem \eqref{sub:1}, which automatically satisfies the inexact optimality condition \eqref{eq:inx}. 
\end{remark}


\subsubsection{Criticality measure convergence rate}
In the previous section, we have proved that the minimum norm subgradients  do not converge to zero when evaluated at $x_k$, but rather when evaluated at a sequence $y_k$  close to $x_k$ (hence, $y_k$, not $x_k$, converges to stationary points). \textit{In this section, we show that the  criticality measure introduced e.g., in \cite{CaGoTo:13,CaGoTo:22,Yua:85} and evaluated at $x_k$ converges to zero and derive explicit rate (consequently, we show that $x_k$ converges to weak stationary points)}.  Let us first introduce the following criticality measure for problem  \eqref{ch:eq:problem}, see also \cite{CaGoTo:13,CaGoTo:22,Yua:85}:

\vspace{-0.5cm}

\begin{align*}
 &\mathcal{M}_f^{r}(x): = f(x) - \min_{\|y - x\|\leq r}\Big( g\big(F(x) + \nabla F(x)(y - x)\big) + h(y) \Big).
\end{align*}

\vspace{-0.2cm}

\noindent First, let us recall the following lemma that connects this criticality measure and the set of weak stationary points (see  Lemma 2.1 in \cite{Yua:85}).
\begin{lemma}\label{lem:yan}
For any \(r \!>\! 0\), the criticality measure \(\mathcal{M}_f^{r}(x^*) \!=\! 0\) if and only if $Df(x^*;d) \geq 0$ for all $d \!\in\! \mathbb{R}^n$ and  for $r \leq \Bar{r} $ we have $\mathcal{M}_f^{r}(x) \leq \mathcal{M}_f^{\bar{r}}(x)$ for all $x\in \emph{dom}\, f$. Additionally, $\mathcal{M}_f^{r}(\cdot)$ is continuous provided that $f$ is continuous and $r \geq 0$. 
\end{lemma}
Clearly, $\mathcal{M}_f^{r}(x^*) \geq 0$ for all $r > 0$. For some $r_k>0$, let us denote:

\vspace{-0.4cm}

\begin{align*}
    \bar{x}_{k+1} &= \argmin_{\|y - x_{k+1}\|\leq r_k} g\Big(F(x_{k+1}) + \nabla F(x_{k+1})(y - x_{k+1})\Big) + h(y).
\end{align*}
\begin{align}\label{eq:meas_subp}
&\mathcal{M}_{s_M(\cdot;x_k)}^{r_k}(x_{k+1}) = g\left(T_p^F(x_{k+1};x_k)\right) + h(x_{k+1}) \\ \nonumber
 &\qquad\qquad\qquad\quad - \min_{\|y - x_{k+1}\|\leq r_k}\Big( g\left(T_p^F(x_{k+1};x_k) + \nabla T_p^F(x_{k+1};x_k)(y - x_{k+1})\right)\\ \nonumber
 &\quad\qquad\qquad\qquad\qquad   +\frac{M}{p!} \|x_{k+1} - x_k\|^{p-1}(x_{k+1} - x_k)^{T}(y - x_{k+1})  + h(y) \Big).
\end{align}


\noindent  Note that if $x_{k+1}$ is a local minimum (or weak stationary point) of subproblem \eqref{sub:1} (i.e.,  $\min_y s_M(y;x_k)$), then according to Lemma \ref{lem:yan}, we have $\mathcal{M}_{s_M(\cdot;x_k)}^{r_k}(x_{k+1}) = 0$ for any $r_k > 0$. Next lemma provides another situation of finding $(r_k,x_{k+1})$, with $x_{k+1}$ not necessarily local minimum/weak stationary point, having  $\mathcal{M}_{s_M(\cdot;x_k)}^{r_k}(x_{k+1})$   small.
\begin{lemma}
Let $x_{k+1}^*$ be a global minimum of  \eqref{sub:1} such that $x_{k+1}^* \not= x_k$ and, additionally, assume that each derivative of order $j=1:p$ of $F_i$ with $i=1:m$, $\nabla^j F_i$, is Lipschitz continuous with constant $L_j^i$ and $h$ is continuous. \red{Then,   there exist $x_{k+1}$, $r_k>0$  and $\delta_k > 0$ satisfying}: 
\begin{align}\label{eq:inx_cr}
\mathcal{M}_{s(\cdot;x_k)}^{r_k}(x_{k+1})\leq \frac{\delta_k}{(p+1)!}\norm{x_{k+1} - x_k}^p  \;\;\; \forall k\geq 0.
\end{align}
\end{lemma}

\vspace{-0.3cm}

\begin{proof}
    From the definition of $x_{k+1}^*$ and Taylor's theorem, we get for all $d\in\mathbb{R}^n$:
    \begin{align*}
        0& \leq  s_M(x_{k+1}^* + d;x_k) - s_M(x_{k+1}^*;x_k)\\
         & = g\left(\sum_{j=0}^{p} \nabla^j T_p^F (x_{k+1}^*;x_k)[d]^j \right) + \sum_{j=0}^{p}\frac{M}{(p+1)!} \nabla^j(\|x_{k+1}^* - x_k\|^{p+1})[d]^j\\
         & \qquad +  \frac{M}{(p+1)!}\nabla^{p+1}(\|x_{k+1}^* - x_k + \tau d\|^{p+1})[d]^{p+1} + h(x_{k+1}^* + d) \\
         &\qquad - g(T_p^F(x_{k+1}^*;x_k)) - \frac{M}{(p+1)!}\|x_{k+1}^* - x_k\|^{p+1} - h(x_{k+1}^*) \\
         & \leq g\left(T_p^F(x_{k+1}^*;x_k) + \nabla T_p^F(x_{k+1}^*;x_k)[d] \right) + g\left(\sum_{j=0}^{p} \nabla^j T_p^F (x_{k+1}^*;x_k)[d]^j\right)\\
         & \qquad + \sum_{j=1}^{p}\frac{M}{(p+1)!} \nabla^j(\|x_{k+1}^* - x_k\|^{p+1})[d]^j + \frac{M}{(p+1)!}\|d\|^{p+1}+ h(x_{k+1}^* + d)\\
         &\qquad - g(T_p^F(x_{k+1}^*;x_k)) - h(x_{k+1}^*) - g\left(T_p^F(x_{k+1}^*;x_k) + \nabla T_p^F(x_{k+1}^*;x_k)[d] \right).
    \end{align*}
for some scalar  $\tau \in (0, \,1)$.    Since $g$ is $L_g$-Lipschitz, we further get:
    \begin{align*}
      &g(T_p^F(x_{k+1}^*;x_k)) + h(x_{k+1}^*) - g\left(T_p^F(x_{k+1}^*;x_k) + \nabla T_p^F(x_{k+1}^*;x_k)[d] \right) - h(x_{k+1}^* + d)\\
      & \leq \! \sum_{j=1}^{p} \! \frac{M}{(p\!+\!1)!} \nabla^j(\|x_{k+1}^* \!-\! x_k\|^{p+1})[d]^j \!+\!  \frac{M}{(p\!+\!1)!}\|d\|^{p+1}  \!+\! L_g  \| \sum_{j=2}^{p} \! \nabla^j T_p^F (x_{k+1}^*;x_k)[d]^j\|.
    \end{align*}
 Since $x_{k+1}^* \not= x_k$, then  $r_k = \|x_{k+1}^* - x_k\| >0$. Then, for any  $\|d\| \leq r_k$,  we obtain:
    \begin{align*}
 &g(T_p^F(x_{k+1}^*;x_k)) + h(x_{k+1}^*) - g\left(T_p^F(x_{k+1}^*;x_k) + \nabla T_p^F(x_{k+1}^*;x_k)[d] \right) - h(x_{k+1}^* + d)\\
      &\leq \sum_{j=1}^{p} M \|x_{k+1}^* - x_k\|^{p+1-j}\|d\|^j + \frac{M}{(p+1)!}\|d\|^{p+1} +L_g \left\| {\sum_{j=2}^{p} \nabla^j T_p^F (x_{k+1}^*;x_k)[d]^j}\right \| \\
      &\leq Mp\|x_{k+1}^* - x_k\|^{p+1} + \frac{M}{(p+1)!}\|x_{k+1}^* - x_k\|^{p+1} + L_g\sum_{j=2}^{p}\|\nabla^j T_p^F (x_{k+1}^*;x_k)\| \|d\|^j\\
      &\leq \Big( \frac{M}{(p+1)!}\|x_{k+1}^* - x_k\| +  L_g\sum_{j=2}^{p}\|\nabla^j T_p^F (x_{k+1}^*;x_k)\| \|x_{k+1}^* - x_k\|^{j-p} \\
      & \qquad +  Mp\|x_{k+1}^* - x_k\| \Big)  \cdot\|x_{k+1}^* - x_k\|^p := \delta_k\|x_{k+1}^* - x_k\|^p,
    \end{align*}
     where we used  that $\|\nabla^j T_p^F (x_{k+1}^*;x_k)\| \leq \sum_{l=j}^{p}\frac{1}{(l-j)!}\|\nabla^l F(x_k)\|  \|x_{k+1}^* - x_k\|^{l-j} \leq \max\limits_{j=2:p} \sqrt{\sum_{i=1}^{m} (L_{j-1}^i)^2 }  \sum_{l=j}^{p}\frac{1}{(l-j)!}   \|x_{k+1}^* - x_k\|^{l-j}$. Defining  $y= x_{k+1}^* + d$  and using the definition of the criticality measure  $\mathcal{M}_{s_M(\cdot;x_k)}^{r_k}$  given in \eqref{eq:meas_subp}, we obtain:
    \begin{align*}
    \mathcal{M}_{s_M(\cdot;x_k)}^{r_k}(x_{k+1}^*)\leq  \delta_k\|x_{k+1}^* - x_k\|^p.  
    \end{align*}
Finally,  continuity of $\mathcal{M}_{s_M(\cdot;x_k)}^{r_k}(\cdot)$ ensures the existence of a feasible neighborhood of $x_{k+1}^*$ in which $x_{k+1}$ can be chosen satisfying \eqref{eq:inx_cr}. This completes our proof.
\end{proof}

\noindent The next theorem derives convergence rate for the iterates $x_k$ of RHOTA in the criticality measure $\mathcal{M}_f^r(x_k)$.
\begin{theorem}
\label{th4.6}
Let the assumptions of Theorem \ref{th:1} hold and $(x_k)_{k\geq 0}$ be generated by RHOTA such that $x_{k+1}$ satisfies the  inexact optimality condition \eqref{eq:inx_cr} for subproblem \eqref{sub:1} for  given $0 < r_{\min} \leq r_{k} \leq r_{\max}$ and $ 0 \leq \delta_k \leq \delta_{\max}$. Denote $\bar{C}\!=\! \frac{(p+1)!(f(x_0) - f^*)}{M - L_g\|L_p\|}$ and  $\bar{D} = \frac{\left(2\bar{C}^{\frac{1}{p+1}} + (p+1)r_{\max}\right)L_g\|L_p\| + (p+1)M r_{\max} + (p+1)\delta_{\max}}{(p+1)!}$. Then, the following convergence rate holds:
\begin{align*}
 \min_{j=0:k} \mathcal{M}_f^{r_{\min}}(x_j) \leq \frac{\bar{D}\Big((f(x_0) - f^*)(p+1)!\Big)^{\frac{p}{p+1}}}{(M - L_g\|L_{p}\|)^{\frac{p}{p+1}}(k+1)^{\frac{p}{p+1}}}  , 
\end{align*}
\end{theorem}
\begin{proof}
From the definitions of the criticality measure $\mathcal{M}_f$ and of $\bar{x}_{k+1}$, we have:
\begin{align*}
 \mathcal{M}_f^{r}&(x_{k+1}) = g(F(x_{k+1})) - g\big(F(x_{k+1}) + \nabla F(x_{k+1})(\bar{x}_{k+1} - x_{k+1})\big)- h(\bar{x}_{k+1})\\
    & = g(F(x_{k+1})) - g\left((T_p^F (x_{k+1};x_k)\right) - g\left(F(x_{k+1}) + \nabla F(x_{k+1})(\bar{x}_{k+1} - x_{k+1})\right) \\
   &\qquad + g\left(T_p^F(x_{k+1};x_k) + \nabla T_p^F(x_{k+1};x_k)[\bar{x}_{k+1} - x_{k+1}]\right) + g\left(T_p^F (x_{k+1};x_k)\right)\\
     &\qquad - g\left(T_p^F(x_{k+1};x_k) + \nabla T_p^F(x_{k+1};x_k)[\bar{x}_{k+1} - x_{k+1}]\right)\\
       &\qquad - \frac{M}{p!}\|x_{k+1} - x_k\|^{p-1}(x_{k+1} - x_k)^{T}(\bar{x}_{k+1} - x_{k+1})\\
       &\qquad +  \frac{M}{p!}\|x_{k+1} - x_k\|^{p-1}(x_{k+1} - x_k)^{T} (\bar{x}_{k+1} - x_{k+1}) - h(\bar{x}_{k+1})\\
       &\leq \frac{L_g\|L_p\|}{(p+1)!}\|x_{k+1} - x_k\|^{p+1} + L_g\|F(x_{k+1}) + \nabla F(x_{k+1})(\bar{x}_{k+1} - x_{k+1}) \\
       &\qquad -  T_p^F(x_{k+1};x_k) - \nabla T_p^F(x_{k+1};x_k)[\bar{x}_{k+1} - x_{k+1}]\|\\
       &\qquad + \frac{M}{p!}\|x_{k+1} - x_k\|^p\|\bar{x}_{k+1} - x_{k+1}\| + \mathcal{M}_{s(\cdot;x_k)}^{r_k}(x_{k+1})\\
       &\leq \frac{2 L_g\|L_p\|}{(p+1)!}\|x_{k+1} - x_k\|^{p+1} + \frac{L_g\|L_p\| + M}{p!}\|x_{k+1} - x_k\|^{p}\|\bar{x}_{k+1} - x_{k+1}\|\\
       &\qquad +\mathcal{M}_{s(\cdot;x_k)}^{r_k}(x_{k+1}) \quad  \forall k \geq 0,
\end{align*}
where the first inequality follows from the fact that $g$ is $L_g$ Lipschitz, from inequality \eqref{main:eq}, the Cauchy-Schwartz inequality and from the fact that:
\begin{align*}
   &g(T_p^F(x_{k+1};x_k)) - \Big( g\left(T_p^F(x_{k+1};x_k) + \nabla T_p^F(x_{k+1};x_k)(\bar{x}_{k+1} - x_{k+1})\right)\\
 &   + \frac{M}{p!}\|x_{k+1} - x_k\|^{p-1}(x_{k+1} - x_k)^T(\bar{x}_{k+1} - x_{k+1})  + h(\bar{x}_{k+1}) \Big) \\
 & \leq   g(T_p^F(x_{k+1};x_k))  -\min_{\|y - x_{k+1}\|\leq r_k} \Big( g\left(T_p^F(x_{k+1};x_k) + \nabla T_p^F(x_{k+1};x_k)(y - x_{k+1})\right)\\
 &\quad   \frac{M}{p!}\|x_{k+1} - x_k\|^{p-1}(x_{k+1} - x_k)^T (y - x_{k+1})  + h(y) \Big) = \mathcal{M}_{s(\cdot;x_k)}^{r_k} (x_{k+1}). 
\end{align*}
From the decrease \eqref{eq:dec}, we get:
\begin{align*}
\|x_{k+1} - x_k\|^{p+1}\leq \frac{(p+1)!}{M - L_g\|L_p\|}(f(x_k) - f(x_{k+1}))&\leq \frac{(p+1)!}{M - L_g\|L_p\|}(f(x_0) - f^*)=\bar{C}.    
\end{align*}
Then, from Lemma \ref{lem:yan} we further obtain:
\begin{align*}
\mathcal{M}_f^{r_{\min}}(x_{k+1})&\leq 
\mathcal{M}_f^{r_k}(x_{k+1})\\
&\leq    \frac{(2\bar{C}^{\frac{1}{p+1}} + (p+1)r_k)L_g\|L_p\| + (p+1)Mr_k + (p+1)\delta_k}{(p+1)!}\|x_{k+1} - x_k\|^p \\
&\leq \bar{D}\|x_{k+1} - x_k\|^p.
\end{align*}
Finally, using the definition of $\bar{k}$ given in \eqref{eq:cv_pnt},  we get:
\begin{align*}
  \min_{j=0:k}\mathcal{M}_f^{r_{\min}}(x_j)^{\frac{p+1}{p}}\leq \bar{D}^{\frac{p+1}{p}}\|x_{\bar{k}+1} - x_{\bar{k}}\|^{p+1}  \stackrel{\eqref{eq:cv_pnt}}{\leq} \bar{D}^{\frac{p+1}{p}}\frac{(f(x_0) - f^*)(p+1)!}{(M - L_g\|L_{p}\|)(k+1)},
\end{align*}
which yields our statement.
\end{proof}
\begin{remark}
Note that if \(x_{k+1}\) is a (local) minimum (or just a weak stationary point)  of the subproblem \eqref{sub:1}, then \(\mathcal{M}_{s(\cdot;x_k)}^{r_k}(x_{k+1}) = 0\) and thus condition \eqref{eq:inx_cr} holds. In Section \ref{sec:sol} we show that we can compute  a global minimum of \eqref{sub:1}  within RHOTA  for the particular case  $g(\cdot)=\|\cdot\|$ and $p=2$ (consequently, satisfying \eqref{eq:inx_cr}). Moreover, from Theorem \ref{th4.6} one can see that the original sequence $x_{k}$ generated by RHOTA does converge to weak stationary points of the original problem \eqref{ch:eq:problem}.
\end{remark}


\subsection{Improved convergence rate under KL}
In this section, \textit{we establish improved convergence rates for RHOTA algorithm under the KL property, i.e., we prove linear/sublinear convergence in function values for the original sequence $(x_k)_{k\geq 0}$ generated by RHOTA}. We denote the set of limit points of  $(x_{k})_{k\geq 0}$ by $\Omega (x_{0})$:

\vspace{-0.6cm}

\begin{align*}
\Omega (x_{0}) = &\lbrace \bar{x}\in \mathbb{R}^n: \exists (k_{t})_{t\geq0} \text{ $\nearrow$ },
\text{ such that } x_{k_{t}}\to \bar{x} \text{ as } t\to \infty \rbrace.
\end{align*} 

\vspace{-0.3cm}

\noindent The next lemma derives some properties for $\Omega (x_{0})$.

\begin{lemma}\label{lem:3}
Let the assumptions of Theorem \ref{th:2} hold.  Additionally, assume that $(x_k)_{k\geq 0}$ is bounded and $f$ is continuous. Then, we have: $\emptyset \neq \Omega(x_{0}) \subseteq\texttt{stat}\, f $, $ \Omega(x_{0})$ is compact and connected set, and $f$ is constant on $ \Omega(x_{0})$, i.e., $f(\Omega(x_{0})) = f_*$.
\end{lemma}
\begin{proof}
Let us prove that $f(\Omega(x_{0}))$ is constant. From the descent \eqref{eq:dec} we have that $\left( f(x_{k})\right)_{k\geq 0}$ is monotonically decreasing, and since $f$ is assumed to be bounded from below, it converges. Let us say to $f_*>-\infty $, i.e., $f(x_{k})\to f_*$ as $k\to \infty$. On the other hand, let $x_{*}$ be a limit point of the sequence $\left(x_{k}\right)_{k\geq0}$. This means that there exists a subsequence $\left( x_{k_{t}}\right)_{t\geq0}$ such that $x_{k_{t}}\to x_{*}$. Since $f$ is continuous, we get $f(x_{k_{t}})\to f(x_{*}) = f_*$ and hence, we have $f(\Omega(x_{0}))=f_*$. The closeness property of $\partial f$ implies that $S_{f}(x_{*})=0$, and thus $0\in \partial f(x_{*})$. This proves that $x_{*}$ is a stationary point of $f$ and thus $\Omega(x_{0})$ is nonempty. By observing that $\Omega(x_{0})$ can be viewed as an intersection of compact sets, $\Omega(x_{0})=\cap_{q\geq 0} \overline{\cup_{k\geq q}\lbrace x_{k}\rbrace}$ so it is also compact. The connectedness follows from \cite{BolSab:14}. This completes the proof.  
\end{proof}


\noindent Next, we derive improved convergence rates in function values for the sequence $(x_k)_{k\geq 0}$ generated by RHOTA, not for the artificial sequence $(y_k)_{k\geq 0}$ as in Theorem \ref{th:2}.

\begin{theorem}\label{th:3}
Let the assumptions of Lemma \ref{lem:3} hold. Additionally, assume that $f$ satisfy the KL property \eqref{eq:kl} on $\Omega(x_0)$. Then, the following convergence rates hold for $(x_{k})_{k\geq 0}$ generated by  RHOTA algorithm for $k$ sufficiently large:
\begin{enumerate}
\item If $q\geq\frac{p+1}{p}$, then $f(x_{k})$ converges to ${f_*}$ linearly. 
\item If $q < \frac{p+1}{p} $, then $f(x_{k})$ converges to ${f_*}$ at sublinear rate of order $\mathcal{O}\left(\frac{1}{k^{\frac{pq}{p+1-pq}}}\right)$.
\end{enumerate}  
\end{theorem}
\begin{proof}
We have:
\begin{align*}
    &f(x_{k + 1}) - f_*\stackrel{\eqref{main:eq}}{\leq} s_M(x_{k+1};x_k) + \frac{L_g\|L_p\| - M}{(p+1)!}\|x_{k+1} - x_k\|^{p+1}  - f_* \\
    &\stackrel{\eqref{eq:inx}}{\leq} \min_{ {y: \; \|y - x_k\|\leq D_k }} s_M(y;x_k) + \frac{\delta + L_g\|L_p\| - M}{(p+1)!}\norm{x_{k+1} - x_k}^{p+1} - f_* \\
    &\stackrel{\eqref{main:eq}}{\leq}\!\!\!\! \min_{ y: \; \|y - x_k\|\leq D_k }\!\! f(y) \!-\! f_*\! +\! \frac{M \!+\! L_g\|L_{p}\|}{(p\!+\!1)!} \|y - x_k\|^{p+1} \!+\! \frac{ \delta + L_g\|L_{p}\|\!-\!M}{(p\!+\!1)!}\|x_{k+1} \!-\! x_k\|^{p+1} \\
    &\leq  f(y_{k+1})  - f_*  \!+\! \left(\frac{M + L_g\|L_p\|}{(p+1)!} \|y_{k+1} - x_k\|^{p+1} \!+\! \frac{\delta+ L_g\|L_p\|\!-\!M }{(p+1)!}\|x_{k+1} - x_k\|^{p+1}\right)\\
    &\leq \sigma_q\text{dist}(0,\partial f(y_{k+1}))^q  + \left(\frac{M + L_g\|L_p\| + L_\mu(\delta +L_g\|L_p\|\!-\!M )}{(p+1)!}\|x_{k+1} - x_k\|^{p+1}\right)\\ 
    &\leq  \left(\frac{\sigma_q (L_{\mu})^{\frac{pq}{p+1}} \mu^{q}}{(p!)^{q}} \|x_{k+1} - x_k\|^{pq} \frac{(1\!-\!L_\mu)M \!+\! (1\!+\!L_\mu)L_g\|L_p\| + L_\mu \delta}{(p+1)!}\|x_{k+1} - x_k\|^{p+1}
    \right)\\
    &\stackrel{\eqref{eq:dec}}{\leq} C_1\big(f(x_{k}) - f(x_{k+1})\big)^{\frac{pq}{p+1}} + C_2\big(f(x_{k}) - f(x_{k+1})\big),
\end{align*}
where the fourth inequality follows from $\|y_{k+1} - x_k\|\leq D_k$, the fifth inequality is deduced from \eqref{eq:kl} combined with the first assertion of Theorem \ref{th:2}, (i.e., the sequences $(x_k)_{k\geq 0}$ and $(y_k)_{k\geq 0}$ share the same limit point) and the sixth inequality follows from \eqref{eq:seq-y} combined with the first assertion of Theorem \ref{th:2}. Here $ C_{1}= \frac{\sigma_q (L_{\mu})^{\frac{pq}{p+1}} \mu^{q}}{(p!)^{q}} \left(\frac{(p+1)!}{M - L_g\|L_p\|}\right)^{\frac{pq}{p+1}}$ and $C_2 = \frac{M + L_g\|L_p\| + L_\mu(\delta +L_g\|L_p\|\!-\!M )}{M - L_g\|L_p\|}$. Let us  denote $\Delta_{k} = f(x_{k}) - f_*$. Subsequently, we derive the following recurrence:
\begin{align*}
    \Delta_{k+1} \leq C_{1}\left( \Delta_{k}-\Delta_{k+1}\right)^{\frac{qp}{p+1}} + C_{2}\left( \Delta_{k}-\Delta_{k+1}\right).
\end{align*}
Using Lemma 6 in \cite{NabNec:22} with $\theta = \frac{p+1}{pq}$, our assertions follow. 
\end{proof}
\begin{remark}
In this section, we have derived improved convergence rates in terms of function values for  sequence $(x_k)_{k\geq 0}$ generated by RHOTA,  by leveraging higher-order information to solve  problem  \eqref{ch:eq:problem}, and \textit{to our knowledge, these rates represent novel findings for such problems  when employing higher-order information}. Notably, for $p=1$, our results align with the convergence rates in \cite{NaToNe:23}.
\end{remark}

\subsection{Adaptive regularized higher-order Taylor approximation method}\label{sec:ad}
In RHOTA algorithm, we need to compute a regularization parameter $M > L_g\|L_{p}\|$. However, in practice, determining Lipschitz constants $L_g$ and $L_{p}$ may be challenging. Consequently, in this section, we introduce an adaptive regularized higher-order Taylor algorithm (A-RHOTA), which does not require prior knowledge of these constants.
\begin{algorithm}
\label{alg:1}
\caption{A-RHOTA algorithm}
\begin{algorithmic}
\STATE{Given $x_0$ and $M_0, R_0 >0$  and $i, k=0$.}
\WHILE{some criterion is not satisfied}
\STATE{1. \; define $s_{M_k}(y;x_k) := g\left(T_{p}^F(y;x_k)\right) + \frac{2^{i}M_k}{(p+1)!}\|y - x_k\|^{p+1} + h(y)$.}
\STATE{2. \;  compute $x_{k+1}$ inexact solution of $\min_y s_{M_k}(y;x_k)$.}
\STATE{\textbf{if}}  descent \eqref{eq:dec} \text{holds}, then go to step 3.\\
\STATE{\;\;\; else} set $i=i+1$ and go to step 1.\\
\STATE{\textbf{end if}}
\STATE{3. \; set $k=k+1$, $M_{k+1} = 2^{i-1}M_k$ and $i=0$.}
\ENDWHILE
\end{algorithmic}
\end{algorithm}

\noindent This line search procedure ensures the decrease \eqref{eq:dec} and finishes in a finite number of steps. Indeed, if $M_k \geq R_0 + L_g\|L_{p}\|$, then from inequality \eqref{main:eq}, we get:

\vspace{-0.5cm}

\begin{align*}
g\left(T_{p}^F(x_{k+1};x_k)\right) - f(x_{k+1}) \geq \frac{-L_g\|L_{p}\|}{(p+1)!}\|x_{k+1} - x_k\|^{p+1}.
\end{align*}

\vspace{-0.3cm}

\noindent This implies that:

\vspace{-0.5cm}

\begin{align*}
f(x_k) - f(x_{k+1}) & = s_{M_k}(x_k;x_k) \!-\! f(x_{k+1}) \geq s_{M_k}(x_{k+1};x_k) \!-\! f(x_{k+1})\\
&\geq \frac{M_k \!-\! L_g\|L_{p}\|}{(p+1)!}\|x_{k+1} \!-\! x_k\|^{p+1} \geq \frac{R_0}{(p+1)!}\| x_{k+1} \!-\! x_k\|^{p+1}.
\end{align*}

\vspace{-0.2cm}

\noindent Note also that we have $M_k \leq 2(R_0 + L_g\|L_p\|) $ for all $k\geq 0$. Consequently,  similar arguments as before allows us to derive convergence rates similar to  Theorems \ref{th:2} and \ref{th:3} for  \red{algorithm A-RHOTA based on the inexact optimality condition \eqref{eq:inx}}.


\section{Efficient solution of the nonconvex  subproblem}\label{sec:sol}
In this section, we present an efficient implementation of RHOTA algorithm for the case $g(\cdot) = \|\cdot\|$, $p=2$ and quadratic  $h(x) = (1/2) x^{T}Bx + ax$. Within this context, $x_{k+1}$ is the solution  of the following nonconvex subproblem (see  subproblem  \eqref{sub:1}):
\begin{align}\label{eq:c-sbp}
\mathcal{P}^* = \min_{x \in \mathbb{R}^n}\;   &\left\|F(x_k) + \langle \nabla F(x_k), x - x_k \rangle + \frac{1}{2}\nabla^2F(x_k) [x - x_k]^2\right\| \\
&\quad + \frac{M}{6} \|x - x_k \|^{3}  + \frac{1}{2}(x - x_k)^{T}B(x - x_k) + \langle a + Bx_k, x - x_k \rangle. \nonumber
\end{align}
with $\langle \nabla F(x_k), x - x_k \rangle = \left[\langle \nabla F_1(x_k), x - x_k \rangle,\!\cdots\!,\langle \nabla F_m(x_k), x - x_k \rangle \right]$ and  $\nabla^2 F(x_k) [x - x_k]^2 
 = \left[\nabla^2 F_1(x_k)[x - x_k]^2,\cdots,\nabla^2 F_m(x_k)[x - x_k]^2\right]$. Denote $u=(u_1,\!\cdots\!,u_m)$. Then, this subproblem is equivalent to:
\begin{align*}
     \min_{x \in \mathbb{R}^n} \max_{\substack{\|u\| \leq 1}}\; &\sum_{i=1}^{m} u_i F_i(x_k) + \left\langle \sum_{i=1}^{m}u_i   \nabla F_i(x_k) + a + Bx_k, x - x_k \right\rangle \\
     &+ \frac{1}{2} \left\langle \left(\sum_{i=1}^{m} u_i   \nabla^2 F_i(x_k) + B\right) (x - x_k),(x - x_k) \right\rangle 
    + \frac{M}{6} \|x - x_k\|^{3}.
\end{align*}
Further, the last term can be  written equivalently as:
\begin{align*}
  \frac{M}{6} \|x - x_k\|^{3} = \max_{w\geq 0} \left( \frac{w}{4}\|x - x_k\|^2 - \frac{1}{12 M^2} w^3 \right).
\end{align*}

\noindent Denote $H_k(u,w) = \sum_{i=1}^{m} u_i   \nabla^2 F_i(x_k) + B + \frac{w}{2} I_n$, $G_k(u) = \sum_{i=1}^{m}u_i   \nabla F_i(x_k) + a + Bx_k$ and $l_k(u) = \sum_{i=1}^{m} u_i   F_i(x_k)$. Then, we have:
\begin{align*}
 \mathcal{P}^* =  \min_{x \in \mathbb{R}^n} \max_{ \|u\| \leq 1  w \geq 0}\; &l_k(u) + \left\langle G_k(u), x - x_k\right\rangle + \frac{1}{2} \left\langle H_k(u,w) (x \!-\! x_k), x \!-\! x_k \right\rangle \!-\! \frac{w^3}{12 M^2}.    
\end{align*}
Consider the following notations:
\begin{align*}
&  \theta_k(x,\!u) \!=\! l_k(u) \!+\! \langle G_k(u),x \!\!-\!\! x_k \rangle \!+\! \frac{1}{2}\! \left\langle\!\! \left(\!\sum_{i=1}^{m} \!u_i   \nabla^2 F_i(x_k) \!+\!\! B\!\right)\!\! (x \!-\! x_k),x \!-\! x_k \!\!\right\rangle \!\!+\!\! \frac{M}{6}\|x \!-\! x_k\|^3 \\
& \beta_k(u,w) = l_k(u) -\frac{1}{2} \left\langle H_k(u,w)^{-1} G_k(u) , G_k(u)\right\rangle - \frac{1}{12 M^2} w^3, \;\;\; r_k = \|x_{k+1} - x_k\| \;\; \text{and}\\
    & \mathcal{F}_k = \left\{ (u,w)\in\mathbb{R}^{m}\times \mathbb{R}_+: \|u\| \leq 1\; \text{and}\; \sum_{i=1}^{m}u_i   \nabla^2 F_i(x_k) + B + \frac{w}{2} I\succ 0  \right\}.\\
\end{align*}

\vspace{-0.5cm}

\noindent Then, we have the following theorem:
\begin{theorem}\label{th:4}
If $M >0$, then we have the following relation:
\begin{align*}
    \theta^*:=\min_{x\in\mathbb{R}^n}\max_{\|u\| \leq  1} \theta_k(x,u)  =  \max_{(u,w)\in \mathcal{F}_k}\beta_k(u,w) = \beta^*. 
\end{align*}
For any $(u,w)\in \mathcal{F}_k$  the direction $x_{k+1} = x_k - H_k(u,w)^{-1} G_k(u)$ satisfies:
\begin{align}\label{eq:5.2}
     0\leq \theta_k(x_{k+1},u) - \beta_k(u,w)  = \frac{M}{12} \left(\frac{w}{M} + 2r_k\right)\left(r_k - \frac{w}{M}\right)^2.
 \end{align}
\end{theorem}
\begin{proof}
First, we show  $\theta^*\geq \beta^*$. Indeed, using a similar reasoning as \cite{NesPol:06}, we have:
\begin{align*}
\theta^* & =  \min_{x \in \mathbb{R}^n} \max_{ \substack{(u,w)\in\mathbb{R}^{m}\times \mathbb{R}_{+}\\ \|u\|\leq 1}}\! l_k(u) \!+\! \left\langle G_k(u), x \!-\! x_k \right\rangle + \frac{1}{2} \!\left\langle H_k(u,w)  (x \!-\! x_k), x \!-\! x_k \right\rangle \!- \frac{ w^3}{12 M^2} \\
&\geq  \max_{ \substack{(u,w)\in\mathbb{R}^{m}\!\times  \mathbb{R}_{+}\\ \|u\|\leq 1}} \min_{x \in \mathbb{R}^n} \! l_k(u) \!+\! \left\langle G_k(u), x \!-\! x_k \right\rangle + \frac{1}{2}\! \left\langle H_k(u,w) (x \!-\! x_k),x \!-\! x_k \right\rangle \!- \frac{ w^3}{12M^2} \\
&\geq \max_{(u,w)\in \mathcal{F}_k} \min_{x \in \mathbb{R}^n} l_k(u) \!+\! \left\langle G_k(u), x \!-\! x_k \right\rangle + \frac{1}{2} \left\langle H_k(u,w) (x \!-\! x_k), x \!-\! x_k \right\rangle - \frac{ w^3}{12M^2}\\
& = \max_{(u,w)\in \mathcal{F}_k} l_k(u) -\frac{1}{2} \left\langle H_k(u,w)^{-1} G_k(u) , G_k(u)\right\rangle - \frac{1}{12 M^2} w^3
= \beta^*.
\end{align*}
Further, let $(u,w) \in \mathcal{F}_k$. Then, we have $G_k(u)= - H_k(u,w) (x_{k+1}- x_k)$ and thus:
\begin{align*}
&\theta(x_{k+1},u)  
 = l_k(u)  - \left\langle H_k(u,\!w)(x_{k+1} - x_k), x_{k+1} - x_k \right\rangle \\
&\qquad + \frac{1}{2}\! \left\langle\! \left(\sum_{i=0}^{m} \!u_i \nabla^2 F_i(x_k) + B\right) (x_{k+1} - x_k), x_{k+1} - x_k \right\rangle \!+\! \frac{M}{6}r_k^3  \\
& = l_k(u) \!-\! \frac{1}{2} \left\langle  \left(\sum_{i=0}^{m} u_i \nabla^2 F_i(x_k) + B  + \frac{w}{2} I_n \right)(x_{k+1} - x_k), x_{k+1} - x_k  \right\rangle - \frac{w}{4}r_k^2 +  \frac{M}{6}r_k^3\\
&=\beta_k(u,w) + \frac{1}{12M^2}w^3 - \frac{w}{4}r_k^2 +  \frac{M}{6}r_k^3 =  \beta_k(u,w) + \frac{M}{12}\left(\frac{w}{M}\right)^3 \!\!-\! \frac{M}{4} \left(\frac{w}{M}\right)r_k^2+ \frac{M}{6}r_k^3\\
&= \beta_k(u,w) + \frac{M}{12} \left(\frac{w}{M} + 2r_k\right)\left(r_k \!-\! \frac{w}{M}\right)^2,
\end{align*}
which proves \eqref{eq:5.2}. Note that we have:
\begin{align*}
    \nabla_w \beta_k(u,w) & = \frac{1}{4} \|x_{k+1} - x_k\|^2 - \frac{1}{4M^2} w^2  = \frac{1}{4} \left(r_k + \frac{w}{M}\right)\left(r_k - \frac{w}{M}\right).
\end{align*}
Therefore if $\beta^*$ is attained at some $(u^*,w^*) >0$ from $\mathcal{F}_k$, then $\nabla \beta_k(u^*,w^*) = 0$. This implies $\frac{w^*}{M} = r_k(u^*,w^*)$ and by \eqref{eq:5.2} we conclude that $\theta^* = \beta^*$. 
\end{proof}
\begin{remark}
\textit{The global minimum} of the nonsmooth nonconvex problem \eqref{eq:c-sbp} is:
\begin{align*}
    x_{k+1} = x_k - H_k(u_k,w_k)^{-1}G_k(u_k),
\end{align*}
with $(u_k,w_k)$  solution of the following convex dual problem:
\begin{align}\label{eq:d-sbp}
&\max_{(u,w)\in \mathcal{F}_k}\; l_k(u) - \frac{1}{2} \left\langle H_k(u,w)^{-1} G_k(u) , G_k(u)\right\rangle - \frac{1}{12 M^2} w^3,
\end{align}
i.e., a maximization of a concave function over a convex set $\mathcal{F}_k$. Hence, if $m$ is of moderate size, this \textit{convex dual problem} \eqref{eq:d-sbp} can be solved very efficiently by interior point methods \cite{NesNem:94}. Thus,  RHOTA algorithm is implementable for $p=2$, since we can effectively compute in fact the \textit{global minimum} $x_{k+1}$ of the subproblem \eqref{sub:1} for $g(\cdot)=\|\cdot\|$ using the powerful tools from convex optimization.
\end{remark}


\section{Applications of RHOTA to nonlinear programming}
In this section, we investigate the behavior of RHOTA algorithm when applied to specific  problems of the form \eqref{ch:eq:problem}. First, we consider solving systems of nonlinear equations and then we consider optimization problems with nonlinear equality constraints. \textit{In both cases, we derive new convergence rates that have not been yet considered in the literature.}  

\subsection{Nonlinear least-squares}
\label{sect:nls}
In this section  we focus on the task of solving a system of nonlinear equations (we assume that this system admits solutions):
\begin{align*}
\text{Find}\;x\in\mathbb{R}^n\;\;\text{such that}:\; F_i(x) = 0\;\;\; \forall i=1:m \;\;\; (\text{with} \;  m\leq n). 
\end{align*}
This problem can be reformulated as the following nonlinear least-squares problem:
\begin{align}
\label{eq:sys}
    \min_{x\in\mathbb{R}^n} \|F(x)\|,
\end{align}
which, essentially, represents a particular form of problem  \eqref{ch:eq:problem} with $h=0$ and $g(\cdot)= \|\cdot \|$ (hence, $L_g=1$). For this particular problem, we assume that the statements from Assumption \ref{as:1} related to $F_i$'s hold and additionally,  there exists $\sigma>0$ such that $\sigma_{\min}(\nabla F(x)) \geq \sigma >0$ for all $x \in \mathcal{L}(x_0)$. This nondegenerate  condition has been considered frequently in the literature, see e.g.,  \cite{Nes:07,XiWr:21}. It  holds e.g.,  when the Mangasarian-Fromovitz constraint qualifications (MFCQ) is satisfied (i.e., the gradients of $F_i(x)$'s, $i=1:m$, are linearly independent  for all $x \in \mathcal{L}(x_0)$) and  $\mathcal{L}(x_0)$ is bounded, see \cite{XiWr:21}. Note that to ensure just local convergence for an algorithm around a local minimum  $x^*$, a nondegeneracy  only needs to hold locally around such a solution \cite{NoWr:06}. However, if one wants to prove global convergence for an algorithm, a nondegenerate condition must hold on a set where all the iterates  lie, e.g., on a level set \cite{Nes:07,XiWr:21}. For simplicity, denote $C_\sigma = \frac{M - \|L_p\|}{(p+1)!L_\mu}\left(\frac{p!\sigma}{\mu}\right)^{\frac{p+1}{p}}$. Under this additional nondegeneracy condition, we can establish finite convergence for RHOTA. 


\begin{theorem}
Let  the assumptions of Theorem \ref{th:2} hold for problem \eqref{eq:sys}. Let also $(x_k)_{k\geq 0}$ be generated by RHOTA algorithm and  $(y_k)_{k\geq 0}$ as defined in \eqref{eq:tk1}, and, additionally, assume that $\sigma_{\min}(\nabla F(x)) \geq \sigma >0$ for all $x\in\mathcal{L}(x_0)$. Then,  there exists finite $k \in \{0, 1,\cdots, \bar k\}$, with  $\bar k =\left\lceil \frac{f(x_0)}{C_{\sigma}} \right\rceil $,  such that either $F(x_k) =0$ or $F(y_k) =0$.
\end{theorem}
\begin{proof}
Combining \eqref{eq:dec} with the first statement of Theorem \ref{th:2}, we obtain:
\begin{align*}
f(x_{k+1}) &\leq f(x_k) - \frac{M - \|L_p\|}{(p+1)!}\|x_{k+1} - x_k\|^{p+1} \leq f(x_k) - \frac{M - \|L_p\|}{(p+1)!L_\mu}\|y_{k+1} - x_k\|^{p+1}.
\end{align*}
From the optimality condition of  $y_{k+1}$ (see \eqref{eq:tk1}), we get:
\begin{align*}
-\frac{\mu}{p!}\|y_{k+1} - x_k\|^{p-1}\left(y_{k+1} - x_k\right) \in \partial f(y_{k+1})  =  \nabla F(y_{k+1}) d_{k+1},   
\end{align*}

\vspace{-0.5cm}

\noindent where

\vspace{-0.3cm}

\begin{equation}\label{eq:def_sub}
d_{k+1} \in \partial \|F(y_{k+1})\|=
\begin{cases}
\frac{F(y_{k+1 })}{\|F(y_{k+1 })\|}\;\;\;\;\;\; &\text{if}\;\;\; F(y_{k+1 }) \neq 0\\
\left\{d \in\mathbb{R}^m : \| d \| \leq 1\right\} \;\;\;\;\;&\text{if}\;\;\; F(y_{k+1}) = 0.
\end{cases}
\end{equation}
 We distinguish  two cases. First, given $\bar{k} \geq 1$, consider that for all  $k \in \{0, 1,\cdots, \bar k-1\}$,  $F(y_{k+1}) \neq 0$. Since $(x_k)_{k\geq 0}, (y_k)_{k\geq 0} \subset \mathcal{L}(x_0)$, from the  nondegeneracy condition of the Jacobians (i.e.,  $\| \nabla F(x) d \| \geq \sigma \|d\|$ for any $d \in \mathbb{R}^m$ and $x \in \mathcal{L}(x_0)$), we have:
\begin{align*}
    \frac{\mu}{p!}\|y_{k+1} - x_k\|^p = \frac{\|\nabla F(y_{k+1}) F(y_{k+1})\|}{\|F(y_{k+1})\|} \geq \sigma \quad \forall  k = 0:\bar{k}-1.
\end{align*}
Hence, we obtain constant decrease in function values for the RHOTA iterates:
\begin{align*}
    f(x_{k+1}) \leq f(x_k) - \frac{M - \|L_p\|}{(p+1)!L_\mu}\left(\frac{p!\sigma}{\mu}\right)^{\frac{p+1}{p}}= f(x_k) - C_{\sigma} \quad \forall  k = 0:\bar{k}-1.
\end{align*}
Summing up the last inequality from $0$ to ${k}$, we get:
\begin{align*}
    0\leq f(x_k) \leq f(x_0) - k C_\sigma \quad \forall  k = 0:\bar{k}.
\end{align*}
Thus, if $\bar k = \left\lceil \frac{f(x_0)}{C_{\sigma}} \right\rceil $, we deduce that $0\leq f(x_{\bar k}) \leq f(x_0) - \bar k C_\sigma \leq 0$, or equivalently   $F(x_{\bar k}) = 0$. In the second case  there exists some  $k\in \{0, 1,\cdots, \bar k-1\}$ such that $F(y_{k+1}) = 0$. Together, both cases prove our statement.
\end{proof}

\subsection{Phase retrieval problem}
In this section, we delve into a specific system of nonlinear equations, i.e., we  focus on the task of solving a quadratic system of equations of the form:
\begin{align}\label{eq:sys_qua}
\text{Find}\;x\in\mathbb{R}^n\;\;\text{such that}:\; x^{T}Q_ix = z_i,\; \text{for}\; i=1:m,
\end{align}
where $(Q_i)_{i=1}^m$ are given matrices and $(z_i)_{i=1}^m$ are measurements. Such problems naturally arise in various real-world scenarios, one notable example being phase retrieval problems \cite{CaLiSo:15}. In phase retrieval, optical sensors impose limitations; when illuminating an object $x$, the resulting diffraction pattern is denoted as $\langle a_i,x\rangle$ for $i=1:m$. However, sensors can only measure the magnitude $z_i: = |\langle a_i,x\rangle|^2$ for $i=1:m$. The objective of phase retrieval is to reconstruct the original signal from its magnitude measurements, which can be mathematically stated as \cite{CaLiSo:15}:
\begin{align}
\text{Find}\; x\in\mathbb{R}^n \; (\text{or}\;\mathbb{C}^n)\;\;\; \text{s.t.}: \;\;  z_i = |\langle a_i , x \rangle|^2,\; i=1:m,
\end{align}
where $a_i \in \mathbb{R}^n \; (\text{or} \; \mathbb{C}^n)$ represents a known measurement vector and $z_i$ denotes a known magnitude, for $i=1:m$. When $a,x \in \mathbb{C}^n$, $\langle a , x \rangle := x^{H} a$, with $x^{H}$  the conjugate transpose of $x$. Various approaches have been explored to tackle phase retrieval, with recent research focusing on non-convex methods. For instance, in \cite{CaLiSo:15}, the authors introduce the Wirtinger flow, a gradient-based method that performs descent on the objective function: 

\vspace{-0.5cm}

\begin{align*}
   x\mapsto\frac{1}{2m}\sum_{i=1}^{m} \Big(|\langle a_i , x \rangle|^2 - z_i  \Big)^2.  
\end{align*}   

\vspace{-0.3cm}

\noindent Similarly, \cite{WaGiEl:17} proposes a modified objective and applies a descent method to:

\vspace{-0.4cm}

\begin{align*}
    x\mapsto \frac{1}{2m}\sum_{i=1}^{m}  \Big(|\langle a_i, x \rangle| - \sqrt{z_i}  \Big)^2.
\end{align*}

\vspace{-0.3cm}

\noindent Both approaches in \cite{CaLiSo:15,WaGiEl:17} rigorously demonstrate the exact retrieval of phase information from a nearly minimal number of random measurements, achieved through careful initialization using spectral methods. In a recent study \cite{DuRu:19}, the authors address phase retrieval using a nonsmooth $l_1$ norm formulation:

\vspace{-0.4cm}

\begin{align*}
       x\mapsto \frac{1}{2m}\sum_{i=1}^{m} \Big |  |\langle a_i , x \rangle  |^2 - z_i  \Big | = \frac{1}{2m}\left \| x^{H} Q x - z \right \|_1,
\end{align*}

\vspace{-0.2cm}

\noindent where $x^{H}Q x = (x^{H}Q_1 x, \ldots, x^{H}Q_m x)^T$ and $Q_i = a_i a_i^{H}$ for $i=1:m$. This problem formulation can be expressed as a composition $f(x)=g(F(x))$. To address this nonsmooth minimization problem, \cite{DuRu:19} proposes a prox-linear method (equivalent to the RHOTA algorithm presented in this paper for $p=1$). The signal recovery in their procedure requires a stability condition, which is typically satisfied with high probability for suitable designs $Q$, a bound on the operator $\|Q\|^2$, and a well-initialized iterative process. The prox-linear method exhibits quadratic convergence and achieves exact signal recovery if the number of measurements satisfies $m = 2n$. In the following we consider $a_i, x$ in $\mathbb{R}^n$. In this paper, inspired by \cite{DuRu:19}, we consider the following nonsmooth composite minimization problem for solving problem \eqref{eq:sys_qua}:

\vspace{-0.6cm}

\begin{align}\label{pr:ph_ret}
    \min_{x\in\mathbb{R}^n} f(x):= \| x^{T} Q x - z\|,
    \end{align}

\vspace{-0.3cm}

\noindent where $x^{T} Q x := \left(x^{T} Q_1 x,\cdots,x^{T} Q_m x\right)^T$ and $Q_i = a_i a_i^{T}$ for $i=1:m$. In the sequel, we present a higher-order proximal point algorithm (called HOPP) for solving this type of problems and then proceed to analyze its convergence rate and efficiency.

\vspace{-0.2cm}

\begin{algorithm}\label{alg:3}
\caption{HOPP Algorithm}
\begin{algorithmic}
\STATE{Given $x_0$, positive integer $p$ and $ M>0$. For $k\geq 0$ do:}
\STATE{Compute $x_{k+1}$ solution of the following subproblem:}
\begin{align}\label{sub:prox}
x_{k+1} \in \arg\min_{x\in\mathbb{R}^{n}} \left\| x^{T} Q x  - z \right\| + \frac{M}{p+1}\| x - x_k\|^{p+1}.
\end{align}
\end{algorithmic}
\vspace{-0.3cm}
\end{algorithm}

\vspace{-0.3cm}

\noindent  The "argmin" in \eqref{sub:prox} refers to the set of global minimizers. Higher-order proximal point algorithms have been considered recently in the literature. Indeed, the convergence rates have been extensively studied, with works such as \cite{Nes:23} focusing on the convex case and \cite{NecDan:20} investigating the nonconvex scenarios. Notably, the objective in \eqref{pr:ph_ret} is weakly convex with $L_f := 2\|(\|Q_1\|,\cdots,\|Q_m\|)\|$, as established in \cite{DruPaq:19}. Therefore, for $p=1$, the subproblem \eqref{sub:prox} becomes strongly convex when $M > L_f$.  However, if the constant $M <  L_f$, one cannot ensure the convexity of the subproblem  \eqref{sub:prox} for $p=1$. In the sequel, we show that when $L_f$ is difficult to compute, one can still employ convex optimization tools to solve the subproblem  \eqref{sub:prox}  for any $M>0$ and $p=1$ or $p=2$. Indeed, following the same reasoning as in Section \ref{sec:sol}, the \textit{global solution} of the (non)convex proximal subproblem \eqref{sub:prox} for $p=1$ is:
\begin{align*}
    x_{k+1} = x_k - H_{k,1}(u)^{-1} g_k(u),
\end{align*}
where  we denote $H_{k,1}(u) := \sum_{i=1}^{m} 2u_i Q_i + \frac{M}{2} I_n$, $g_k(u) := \sum_{i=1}^{m} 2 u_i Q_i x_k$ and $l_k(u) := \sum_{i=1}^{m} u_i (x_k^{T}Q_i x_k - z_i)$, with $u$ representing the solution of the convex dual problem:
\begin{align*}
&\max_{u\in \mathcal{F}_1} \; l_k(u) - \frac{1}{2} \left\langle H_{k,1}(u)^{-1} g_k(u) , g_k(u)\right\rangle,
\end{align*}
where $\mathcal{F}_1: = \{ u\in\mathbb{R}^m: \|u\|\leq 1\; \text{and}\;   H_{k,1}(u) \succ 0\}$. Similarly, for $p=2$, we have:
\begin{align*}
    x_{k+1} = x_k - H_{k,2}(u,w)^{-1} g_k(u),
\end{align*}
where $H_{k,2}(u,w) := \sum_{i=1}^{m} 2u_i Q_i + \frac{w}{2} I_n$, with $(u,w)$ is the solution of the following convex dual problem:

\vspace{-0.6cm}

\begin{align*}
\max_{(u,w)\in \mathcal{F}_2} \; l_k(u) - \frac{1}{2} \left\langle H_{k,2}(u,w)^{-1} g_k(u) , g_k(u)\right\rangle - \frac{1}{12 M^2} w^3,
\end{align*}

\vspace{-0.2cm}

\noindent where $\mathcal{F}_2: = \{ (u,w)\in\mathbb{R}^m \times \mathbb{R}_+: \|u\|\leq 1\; \text{and}\;   H_{k,2}(u,w) \succ 0\}$. Hence, our algorithm HOPP can be easily implemented  for any $M>0$ and $p=1, 2$ using standard convex optimization tools  (such as interior point methods \cite{NesNem:94}). Next, we establish the global convergence rate to a stationary point for this algorithm:

\begin{theorem}
Let $f$ be given as in \eqref{pr:ph_ret} and let $(x_{x})_{k\geq 0}$ be generated by  HOPP algorithm. Then, we have:

\vspace{-0.5cm}

\begin{align*}
    &\min_{i=0:k} S_{f}(x_{i}) \leq \left(\frac{(M(p+1)^p)^{\frac{1}{p+1}}(f(x_0) - f^*)}{k^{\frac{p}{p+1}}}\right). 
\end{align*}
\end{theorem}

\vspace{-0.3cm}

\begin{proof}
From the definition of $x_{k+1}$ in \eqref{sub:prox}, we get:
\begin{align*}
 &f(x_{k+1}) + \frac{M}{p+1}\|x_{k+1} - x_k\|^{p+1}\leq f(x_k) \;\; \text{and} \;\;  S_{f}(x_{k+1}) \leq M\|x_{k+1} - x_k\|^p.
\end{align*}
Hence, combining the last two inequalities, we get:
\begin{align*}
S_{f}(x_{k+1})^{\frac{p+1}{p}}\leq M^{\frac{p+1}{p}}\frac{p+1}{M}\big( f(x_k) - f(x_{k+1}) \big).    
\end{align*}
Summing up and taking the minimum, we  get our statement.
\end{proof}

\noindent In order to establish rapid local convergence, we introduce an additional assumption that is related to sharpness or  error bound  for  objective function, see \cite{DuRu:19}.

\begin{assumption}\label{ass:quad}
There exists $\lambda >0$ such that for all $x\in\mathbb{R}^n$ the objective function $f$ defined in \eqref{pr:ph_ret}, having the set of global minima  $X^*$, satisfies:
\vspace{-0.5cm}

 \begin{align*}
     f(x) - f(x^*)\geq \sigma_0 \; \emph{dist}(0,X^*)\;\emph{dist}(x,X^*)  \;\;\; \forall x^* \in X^*, \;\; \text{with} \; \sigma_0 >0.
 \end{align*}

\vspace{-0.2cm}
\end{assumption}
This condition has been proved to hold in the context of  phase retrieval, see \cite{DuRu:19}. For instance, it holds when the matrices $Q_i$'s satisfy the following stability condition \cite{DuRu:19}:
\begin{align*}
    \|(Q_ix)^2 - (Q_iy)^2\|\geq \bar{\sigma}_0 \|x-y\|\;\|x+y\| \;\;\; \forall x,y \in\mathbb{R}^n, \; i=1:m, \;\; \text{with} \; \bar{\sigma}_0>0.
\end{align*}
Next, we derive a fast convergence rate for  HOPP algorithm under sharpness.

\begin{theorem}
Let $f$ be defined as in \eqref{pr:ph_ret} and satisfy Assumption \ref{ass:quad}. Moreover, let $(x_k)_{k\geq 0}$ be generated by  HOPP algorithm. Then, we have:
\begin{align*}
      \frac{\emph{dist}(x_k,X^*)}{\emph{dist}(0,X^*)}\leq \left(\frac{\sigma_0(p+1)}{M\emph{dist}(0,X^*)^{p-1}}\right)^{\frac{1}{p}} \left(\frac{M^{\frac{1}{p}} \;\emph{dist}(x_0,X^*)}{(\sigma_0(p+1) \emph{dist}(0,X^*))^{\frac{1}{p}}} \right)^{(p+1)^k}. 
\end{align*}
\end{theorem}

\begin{proof}
Since $x_{k+1}$ is the global minimum of \eqref{sub:prox}, we have:
\begin{align*}
f(x_{k+1}) &\leq \min_{x\in\mathbb{R}^n} f(x) + \frac{M}{p+1}\|x - x_k\|^{p+1}\leq f(x^*) + \frac{M}{p+1}\|x^* - x_k\|^{p+1}.
\end{align*}
Taking the infimum over $x^*\in X^*$, we further obtain:

\vspace{-0.5cm}

\begin{align*}
f(x_{k+1}) - f(x^*) \leq \frac{M}{p+1}\text{dist}(x_k,X^*)^{p+1}.
\end{align*}

\vspace{-0.2cm}

\noindent Combining this inequality with Assumption \ref{ass:quad}, we get:

\vspace{-0.5cm}

\begin{align*}
    \sigma_0 (p+1)\; \text{dist}(0,X^*)\text{dist}(x_{k+1},X^*)\leq M \text{dist}(x_k,X^*)^{p+1}.
\end{align*}


\noindent Dividing each side by $\text{dist}(0,X^*)^{p+1}$, we get:
\begin{align*}
    \frac{\text{dist}(x_{k+1},X^*)}{ \text{dist}(0,X^*)}\leq \frac{M \text{dist}(0,X^*)^{p-1}}{\sigma_0 (p+1) } \left(\frac{\text{dist}(x_k,X^*)}{ \text{dist}(0,X^*)}\right)^{p+1}. 
\end{align*}
\noindent Unrolling the last recurrence, yields our statement.
\end{proof}

\noindent Note that if $M \text{dist}(x_0,X^*)^p < \sigma_0(p+1) \text{dist}(0,X^*)$, then  faster convergence is guaranteed for HOPP iterates with the increasing value of $p$. Note also that for $p=1$, we recover the quadratic convergence rate obtained in \cite{DuRu:19}. Furthermore, \textit{the flexibility in selecting the parameter $M$ is significant}: given an arbitrary initial point $x_0$ (not necessarily close to $X^*$), an appropriate  choice of $M$ (i.e., sufficiently small) guarantees very fast convergence of HOPP iterates to the global minima of \eqref{pr:ph_ret}.


\subsection{Equality constrained nonlinear problems}
\noindent Let us now consider an optimization problem with nonlinear equality constraints:

\vspace{-0.5cm}

\begin{align}\label{eq:no-li-pr}
\min_{x \in  \mathbb{R}^n }\; h(x)\;\;\;  \text{s.t.:} \;\;\;  F(x) = 0,
\end{align}

\vspace{-0.1cm}

\noindent where $F(x) = (F_1(x),\cdots,F_m(x))$, with $m\leq n$, and $h$ is proper lsc function. For a given positive constant $\rho$, the exact penalty reformulation of \eqref{eq:no-li-pr} is \cite{NoWr:06}:

\vspace{-0.5cm}

\begin{align}\label{eq:expen}
    \min_{x \in  \mathbb{R}^n } f(x): =  h(x) + \rho\|F(x)\|,
\end{align}

\vspace{-0.1cm}

\noindent which fits into the  formulation \eqref{ch:eq:problem} with $g(\cdot) = \rho\| \cdot\|$. It is known  that, under proper  constraint qualification conditions and for sufficiently large $\rho$, any   stationary point $x^*$ of the exact penalty problem \eqref{eq:expen} corresponds to a Karush-Kuhn-Tucker (KKT) point of  constrained problem \eqref{eq:no-li-pr}  (i.e., $\exists \lambda^*$ \red{finite} s.t. $0 \in \partial h(x^*) + \nabla F(x^*)^T \lambda^*$ and $F(x^*) =0$), see e.g.,  \cite{CaGoTo:11}. When the objective function,  $h$, exhibits smoothness,  constraint qualifications are naturally related to the constraints themselves, $F$, see e.g.,  LICQ or MFCQ \cite{NoWr:06}. However, when  $h$ takes on a non-smooth character, a shift occurs, necessitating the introduction of new constraint qualifications. This adjustment becomes imperative because the non-smoothness of the objective function has the potential to significantly impact the behavior and satisfaction of the constraints. Hence, in such scenarios, a good understanding of these new constraint qualifications becomes essential to navigate the complexities inherent in optimizing non-smooth objectives within nonlinear programming. Thus, in this section, for  problem \eqref{eq:no-li-pr} we assume that the statements from Assumption \ref{as:1} related to $F_i$'s and $f$ hold and additionally, there exist $\sigma  > 0$ such that the following \textit{new constraint qualification condition} holds: 

\vspace{-0.4cm}

\begin{align}
\label{newCQ}
\sigma \|\lambda\| \leq \text{dist}\left(- \nabla F(x)^{T}\lambda, \partial^{\infty} h(x) \right)  \;\;\; \forall  x \in \mathcal{L}(x_0) \;\; \text{and} \;\;  \lambda \in \partial \|F(x)\|.
\end{align} 

\vspace{-0cm}

\noindent If $h=0$ or $h$ is locally Lipschitz continuous, then $\partial ^{\infty} h(x) = \{0\}$ (see Theorem 9.13 in \cite{Roc:70}) and thus \eqref{newCQ} reduces to the nondegeneracy condition from Section \ref{sect:nls}: $\sigma_{\min}(\nabla F(x)) \geq \sigma$ for all $x \in \mathcal{L}(x_0)$.  A  constraint qualification condition of the form $\sigma \|\lambda\| \leq \text{dist}\left(- \nabla F(x)^{T}\lambda, \partial h(x) \right)$ for all  $x \in \mathcal{L}(x_0)$  \text{and}   $\lambda \in \partial \|F(x)\|$ has been  adopted when analyzing the convergence of  algorithms for solving  \eqref{eq:no-li-pr}, see e.g.,  \cite{Lu:22,SahEft:19}. However, in  \cite{Lu:22,SahEft:19}, the proposed constraint qualification condition loses coherence  when e.g., the nonsmooth component exhibits (local) Lipschitz continuity, such as  $h(\cdot) = \| \cdot \|_1 \; \text{or} \;  \| \cdot \|_2$, while our \eqref{newCQ} imposes only a condition on $\nabla F$ in this case.


\begin{theorem}
\label{th:h-F}
Let the assumptions of Theorem \ref{th:2} hold, and, additionally, the constraint qualification condition \eqref{newCQ} holds. Let $\bar{\rho}>0$ be fixed sufficiently large and the sequence $(x_k)_{k\geq 0}$ be generated by RHOTA  applied to penalty problem \eqref{eq:expen}, with $\rho \geq \bar{\rho}$, $M \geq 2\rho \|L_p\|$,   and $(y_k)_{k\geq 0}$ be given in \eqref{eq:tk1}, with $\mu = 2(M + \rho\|L_p\|)$. Then, any limit point of  the sequence $(x_{k})_{k\geq 0}$ is a KKT point of \eqref{eq:no-li-pr}. Moreover, the convergence  rate towards a KKT point is of order $\mathcal{O}\left(\rho k^{-\frac{p}{p+1}}\right)$.
\end{theorem}

\vspace{-0.2cm}

\begin{proof}
From the optimality conditions of $y_{k+1}$  applied to $f$ given in \eqref{eq:expen} (see \eqref{eq:tk1}), there exists $\lambda_{k+1}\in\partial \|F(y_{k+1})\|$ such that:

\vspace{-0.4cm}

\begin{align}
\label{eq:opt_sub}
\frac{\mu}{p!}\|y_{k+1} - x_k\|^{p-1}(x_k - y_{k+1}) \in \rho \nabla F(y_{k+1})^{T} \lambda_{k+1} + \partial h(y_{k+1}) \quad \forall k \geq 0.
\end{align}

\vspace{-0.2cm}

\noindent This implies that (for simplicity, we denote $\mathcal{N}_{\text{epi}\,h}^{k+1} = \mathcal{N}_{\text{epi}\,h}\big(y_{k+1},h(y_{k+1})\big)$):

\vspace{-0.4cm}

\begin{align*}
&\text{dist}\Big(\left(-\rho \nabla F(y_{k+1})^{T} \lambda_{k+1} , 0\right),   \mathcal{N}_{\text{epi}\,h}^{k+1}\Big)- \left\|\left(\frac{\mu}{p!}\|y_{k+1} - x_k\|^{p-1}(y_{k+1} - x_k),1\right)\right\| \\
\vspace*{-0.6cm}
& \leq \text{dist}\left(\Big(- \rho \nabla F(y_{k+1})^{T} \lambda_{k+1} - \frac{\mu}{p!}\|y_{k+1} \!-\! x_k\|^{p-1}(y_{k+1} \!-\! x_k ) ,-1\Big), \mathcal{N}_{\text{epi}\,h}^{k+1}\right) = 0.   
\end{align*}

\vspace{-0.2cm}

\noindent On the other hand, from the definition of the horizon subdifferential, we get:

\vspace{-0.4cm}

\begin{align*}
 &\text{dist}\left(-\rho \nabla F(y_{k+1})^{T} \lambda_{k+1}, \partial^{\infty}h(y_{k+1}) \right) =  \text{dist}\Big(\big(-\rho \nabla F(y_{k+1})^{T} \lambda_{k+1},0\big), \mathcal{N}_{\text{epi}\,h}^{k+1}\Big).   
\end{align*}

\vspace{-0.2cm}

\noindent  Therefore, combining the last two inequalities with the constraint qualification condition and using that $\partial^{\infty} h(y_{k+1})$ is a cone, we obtain for any $\rho >0$:

 \vspace{-0.4cm}
 
\begin{align}\label{eq:feas}
& \sigma  \rho\|\lambda_{k+1} \| \leq  \text{dist} \left( - \rho\nabla F(y_{k+1})^{T} \lambda_{k+1}, \partial^{\infty} h(y_{k+1}) \right) \leq  \frac{\mu}{p!}\|y_{k+1} - x_k\|^p + 1.  
\end{align}

\vspace{-0.2cm}

\noindent Or, equivalently, using the definition of $M$ and $\rho \geq \bar{\rho}$, we have:

\vspace{-0.4cm}

\begin{align*}
& \|\lambda_{k+1} \| \leq  \left(  \frac{2 M}{\sigma \bar{\rho} p!} + \frac{2 \|L_p\|}{\sigma  p!}\right) \|y_{k+1} - x_k\|^p + \frac{1}{\sigma \bar{\rho}}.  
\end{align*}

\vspace{-0.2cm}

\noindent Since $\|y_{k+1} - x_k\|\to 0$ as $k\to \infty$ (see Theorems \ref{th:1} and \ref{th:2}), then the previous relation implies that for fixed $\bar{\rho}>0$ sufficiently large (e.g., $\sigma \bar{\rho} > 1$) there exists  positive integer $\bar{k}$ such that: 

\vspace{-0.5cm}

\begin{align*}
 \|\lambda_{k+1}\|  < 1  \Longrightarrow  F(y_{k+1}) \stackrel{\eqref{eq:def_sub}}{=} 0 \quad \forall k \geq  \bar{k},\; \rho \geq \bar{\rho}.  
\end{align*} 

\vspace{-0.1cm}

\noindent Hence, feasibility is achieved after a finite number of iterations. Additionally, it also follows from \eqref{eq:opt_sub} that for any $k \geq 0$ there exists $h_{y_{k+1}} \in \partial h(y_{k+1})$ such that:

\vspace{-0.3cm}

$$\| \nabla F(y_{k+1})^{T} (\rho \lambda_{k+1}) + h_{y_{k+1}} \|  = \frac{\mu}{p!}\|y_{k+1} - x_k\|^p  \to 0 \;\;  \text{as} \;\;  k\to \infty.$$

\vspace{-0.1cm}

\noindent Using the closedness of the graph of $\partial h$  and basic limit rules, we deduce that any limit point of  the sequence $(y_{k})_{k\geq 0}$ is a KKT point of \eqref{eq:no-li-pr}. Since the set of limit points of  $(x_{k})_{k\geq 0}$ coincides with the set of limit points of  $(y_{k})_{k\geq 0}$ (see Theorem \ref{th:2}), the first statement follows. 
Further, from Theorem \ref{th:2}, there exists  $\bar{k} \in\{0,\cdots,k\}$ such that:

\vspace{-0.5cm}

 \begin{align*}
S_f(y_{\bar{k} + 1}) \leq \left( 
 \frac{ \left(\frac{\mu + \delta + \rho\|L_p\| - M}{\mu - (M + \rho\|L_p\|)}\right) \mu^{\frac{p+1}{p}} (p+1)!}{(M - \rho\|L_p\|)(p!)^{\frac{p+1}{p}} (k+1) } (f(x_0) - f^*)    
 \right)^{\frac{p}{p+1}}.
 \end{align*}

 \vspace{-0.2cm}
 
\noindent Since Assumption \ref{as:1}.3 holds, then $f(x) = h(x) + \rho\|F(x)\|  \geq f^*$ and, consequently, we have $f(x_0) - f^* = \mathcal{O}(\rho)$. In addition, since $\delta \ll \rho$, we deduce the following bound:

\vspace{-0.5cm}

\begin{align}
\label{eq:rhorate}
S_f(y_{\bar{k} + 1}) \leq\mathcal{O}\left( \frac{\rho}{k^{\frac{p}{p+1}}}\right).
\end{align}

\vspace{-0.2cm}

\noindent Further, combining \eqref{eq:feas} with  first assertion of Th. \ref{th:2} and with eq. \eqref{eq:cv_pnt}, we get:

\vspace{-0.5cm}

\begin{align*}
    \|\lambda_{\Bar{k}+1}\| \leq \frac{\mu (L_\mu)^{\frac{p}{p+1}}}{p!\sigma\rho}\frac{\left((f(x_0) - f^*)(p+1)!\right)^{\frac{p}{p+1}} }{ ((M - \rho\|L_p\|)(k+1))^{\frac{p}{p+1}}} + \frac{1}{\sigma\rho} = \mathcal{O}\left(\frac{1}{k^{\frac{p}{p+1}}}\right) + \frac{1}{\sigma\rho},
\end{align*}

\vspace{-0.2cm}

\noindent where $\mathcal{O}(\cdot)$ does not depend on $\rho$. Hence, for any given $\epsilon$, with $0 < \epsilon < \frac{1}{2}$, and for any  $\rho > \frac{2}{\sigma}$, if $k\geq \mathcal{O}\left(\rho^{\frac{p+1}{p}}\epsilon^{-\frac{p+1}{p}}\right)$, then  there exists $h_{y_{\bar{k}+1}} \in \partial h(y_{\bar{k}+1})$ such that:

\vspace{-0.5cm}

\begin{align*}
 S_f(y_{\bar{k} + 1}) \!=\! \| \nabla F(y_{\bar{k}+1})^{T} \!(\rho \lambda_{\bar{k}+1}) \!+\! h_{y_{\bar{k}+1}} \| \!\leq\! \epsilon\; \text{and} \;
 \|\lambda_{\bar{k} + 1}\| \!\leq\! \epsilon \!+\! \frac{1}{\sigma \rho} \!<\! 1  \Rightarrow  F(y_{\bar{k}+1}) \! \stackrel{\eqref{eq:def_sub}}{=} \!0,
\end{align*}

\vspace{-0.2cm}

\noindent i.e., $y_{\bar{k} + 1}$ satisfies $\epsilon$-KKT conditions (but exact feasibility), i.e., second  statement.
\end{proof}

\begin{remark}
From  previous proof one notice that  in order to guarantee  feasibility,   $\rho$ needs to be sufficiently large, e.g.,  $\rho > \frac{1}{\sigma}$. On the other hand it is known that in exact penalty methods one needs to choose  $\rho$  larger than the norm of Lagrange multiplier associated to a  KKT point of  \eqref{eq:no-li-pr} \cite{NoWr:06}. Let $(x^*,\lambda^*)$ be a KKT point of  \eqref{eq:no-li-pr},~i.e.:
\begin{align*}
 (-\nabla F(x^*)^{T}\lambda^*,-1)\in \mathcal{N}_{\text{epi}\,h}(x^*,h(x^*)).   
\end{align*}

\vspace{-0.2cm}

\noindent This implies that:

\vspace{-0.5cm}

\begin{align*}
    \sigma\|\lambda^*\|\stackrel{\eqref{newCQ}}{\leq} \text{dist}(-\nabla F(x^*)^{T}\lambda^*,\partial^{\infty}h(x^*))\leq 1  \;\; \Rightarrow  \;\; \rho > \frac{1}{\sigma} \geq \|\lambda^*\|, 
\end{align*}

\vspace{-0.2cm}

\noindent i.e., we have established the connection between our lower bound $1/\sigma$ and the known lower bound from literature $\|\lambda^*\|$ on the exact penalty parameter  $\rho$.  \textit{To the best of our knowledge, Theorem \ref{th:h-F} provides  the first convergence results for a higher-order exact penalty method for solving the equality constrained optimization problem  \eqref{eq:no-li-pr}, i.e. finding a near KKT point}. Specifically, for $p$ very large we get a rate of order $\mathcal{O}(\epsilon^{-1})$, while for $p=1$ our rate  $\mathcal{O}(\epsilon^{-2})$ aligns with that previously obtained in e.g.,~\cite{CaGoTo:11}.
\end{remark}


\section{Numerical experiments}
In this section, we test the performance of the proposed algorithms in solving systems of quadratic equations using real data. We consider two applications: the output feedback control and the phase retrieval  problems. The implementation details are conducted using MATLAB R2020b on a laptop equipped with an i5 CPU operating at 2.1 GHz and 16 GB of RAM.

\subsection{Solving output feedback control problems}
\noindent In this section we evaluate the performance of RHOTA algorithm for solving static output feedback control problem \eqref{fw: eq3} using data from the COMPl$_{e}$ib library available at \cite{Lei:16}. Let us note that $\nabla^{2} F$, where $F$ is given in \eqref{fw: eq3}, is constant, and therefore $\nabla F$ is Lipschitz (i.e., $p=1$).  Hence, our algorithm RHOTA with $p=1$ can be used for solving the output feedback control problem with mathematical guarantees for finding stationary points. We can effectively implement  RHOTA algorithm by utilizing the Fr$\acute{e}$chet differentiable of the matrix function  $F$ \cite{FaBo:21}. Thus, at each iteration of  RHOTA  with $p=1$, the following convex subproblem needs to be solved:

\vspace{-0.5cm}

\begin{align*}
(X_{k+1},K_{k+1},Q_{k+1}) &= \argmin_{X \succ 0, Q \succ 0, K} \Big\|F(X_k, K_k, Q_k) + \nabla F(X_k,K_k,Q_k)[X - X_k,\\
\vspace*{-0.5cm}
&\qquad K - K_k, Q - Q_k] \Big\|_F
+ \frac{M}{2}\Big \|[X - X_k, K - K_k, Q - Q_k]\Big \|^2_F,  
\end{align*}

\vspace*{-0.1cm}

\noindent where the directional derivative  is  $\nabla F(X_k,K_k,Q_k)[X - X_k, K - K_k, Q - Q_k] = \nabla_X F (X_k,K_k,Q_k)[X - X_k]  +  \nabla_K F (X_k,K_k,Q_k)[K - K_k]  +  \nabla_Q F (X_k,K_k,Q_k)[Q - Q_k]$:

\vspace{-0.5cm}

\begin{align*}
  &\nabla_X F (X_k,K_k,Q_k)[X - X_k] = (A + BK_k C)^{T}(X - X_k) + (X - X_k)(A + BK_k C)\\
 & \nabla_K F (X_k,K_k,Q_k)[K - K_k] = (B(K - K_k)C)^{T}X_k + X_k(B(K - K_k)C)\\
& \nabla_Q F (X_k,K_k,Q_k)[Q - Q_k] = Q - Q_k.
\end{align*}

\begin{table}[]
\small
{ \centering
\begin{adjustbox}{width = \columnwidth,center}
\begin{tabular}{|l|llll|llll|}
\hline
\multirow{2}{*}{} & \multicolumn{4}{l|}{RHOTA ($p=1$)}           & \multicolumn{4}{l|}{BMIsolver \cite{DiGuMiDi:11}}           \\ \cline{2-9} 
                  & \multicolumn{1}{l|}{iter} & \multicolumn{1}{l|}{time(s)} & \multicolumn{1}{l|}{K} & MEV & \multicolumn{1}{l|}{iter} & \multicolumn{1}{l|}{time(s)} & \multicolumn{1}{l|}{K} & MEV \\ \hline
ac3     ($n_x = 5$)             & \multicolumn{1}{l|}{\textbf{4}}     & \multicolumn{1}{l|}{\textbf{1.08}}    & \multicolumn{1}{l|}{$\begin{pmatrix} 2.7633 &  -0.4060  & -2.6203 &  -0.0605 \\  -0.1880  &  1.5857  & -3.5001   & 1.8552\end{pmatrix}$}  &     \textbf{-0.89}            & \multicolumn{1}{l|}{29}     & \multicolumn{1}{l|}{18}    & \multicolumn{1}{l|}{$\begin{pmatrix} 2.9051  & -0.4423  & -2.7215  &  0.0038\\ -0.1084  &  1.7357  & -3.3988   & 1.9438\end{pmatrix}$ }  &   -0.85           \\ \hline
ac8     ($n_x = 9$)             & \multicolumn{1}{l|}{\textbf{6}}     & \multicolumn{1}{l|}{\textbf{1.8}}    & \multicolumn{1}{l|}{$\begin{pmatrix} 1.0279 & -0.4365  & -1.1585       0.0085 &   0.46237 \end{pmatrix}$}  &     \textbf{ -0.44}            & \multicolumn{1}{l|}{43}     & \multicolumn{1}{l|}{5.6}    & \multicolumn{1}{l|}{$\begin{pmatrix} 1.0279 & -0.4365  & -1.1585       0.0085 &   0.46237 \end{pmatrix}$ }  & \textbf{ -0.44}            \\ \hline
cm1\_is    ($n_x = 20$)            & \multicolumn{1}{l|}{\textbf{12}}     & \multicolumn{1}{l|}{\textbf{5.7}}    & \multicolumn{1}{l|}{$\begin{pmatrix} -19.98 \\  -10.97\end{pmatrix}$}  &   \textbf{-4.3e$^{-3}$}             & \multicolumn{1}{l|}{22}     & \multicolumn{1}{l|}{10.6}    & \multicolumn{1}{l|}{$\begin{pmatrix} -17.85 \\ -22\end{pmatrix}$}  & -4.2$e^{-3}$              \\ \hline
cm2\_is     ($n_x = 60$)              & \multicolumn{1}{l|}{\textbf{19}}     & \multicolumn{1}{l|}{\textbf{533}}    & \multicolumn{1}{l|}{$\begin{pmatrix} -5 \\ -7.87 \end{pmatrix}$}  &         \textbf{-1.07e$^{-2}$}     & \multicolumn{1}{l|}{ 114}  & \multicolumn{1}{l|}{2691}    & \multicolumn{1}{l|}{$\begin{pmatrix} -5.6 \\ -7.18 \end{pmatrix}$}  &     -1.02e$^{-2}$      \\ \hline
dis1     ($n_x = 8$)           & \multicolumn{1}{l|}{\textbf{24}}     & \multicolumn{1}{l|}{\textbf{7.6}}    & \multicolumn{1}{l|}{$\begin{pmatrix}  3.125    &   2.817   &  -7.584   &   -5.446 \\   2.817  &     3.71   &   -4.244   &    -4.256\\   -7.584   &   -4.244   &  -0.435    &   1.352\\  -5.446  &    -4.256   &   1.352   &    2.877 \end{pmatrix}$}  &        \textbf{-1.363}          & \multicolumn{1}{l|}{100}     & \multicolumn{1}{l|}{13.4}    & \multicolumn{1}{l|}{$\begin{pmatrix}  3.125    &   2.817   &  -7.584   &   -5.446 \\   2.817  &     3.71   &   -4.244   &    -4.256\\   -7.584   &   -4.244   &  -0.435    &   1.352\\  -5.446  &    -4.256   &   1.352   &    2.877 \end{pmatrix}$}  &      \textbf{-1.363}         \\ \hline
dlr2   ($n_x = 40$)                & \multicolumn{1}{l|}{\textbf{7}}     & \multicolumn{1}{l|}{\textbf{21.6}}    & \multicolumn{1}{l|}{$\begin{pmatrix}1.85 &  1.06 \\ -0.27 & -2.09 \end{pmatrix}$}  &   \textbf{-5e$^{-3}$}               & \multicolumn{1}{l|}{120}     & \multicolumn{1}{l|}{477}    & \multicolumn{1}{l|}{$\begin{pmatrix} -5.6 & 5.5\\  5.5 &  -1.4\end{pmatrix}$}  &   \textbf{-5e$^{-3}$}           \\ \hline
eb1      ($n_x = 10$)           & \multicolumn{1}{l|}{\textbf{8}}     & \multicolumn{1}{l|}{\textbf{3.5}}    & \multicolumn{1}{l|}{-0.551}  &   \textbf{-0.066}               & \multicolumn{1}{l|}{26}     & \multicolumn{1}{l|}{9.2}    & \multicolumn{1}{l|}{-47.377}  &       -0.0212      \\ \hline

he1 ($n_x = 4$) & \multicolumn{1}{l|}{\textbf{5}}    & \multicolumn{1}{l|}{\textbf{1.4}}    & \multicolumn{1}{l|}{$\begin{pmatrix}  0.981  \\  4.469\end{pmatrix}$}  &            -0.22     & \multicolumn{1}{l|}{12}     & \multicolumn{1}{l|}{2.1}    & \multicolumn{1}{l|}{ $\begin{pmatrix}  0.883 \\ 4.022 \end{pmatrix}$}  &    -0.22             \\ \hline
he4  ($n_x = 8$)               & \multicolumn{1}{l|}{\textbf{15}}     & \multicolumn{1}{l|}{\textbf{8.5}}    & \multicolumn{1}{l|}{  $\begin{pmatrix} -2.12  &  3.87 &   1.47 &   -0.26 &  -0.04 &   0.65 \\  3.75 &  -14.55 &   -1.48 &    1.14 &    5.35 &    2.03\\ -0.67 &    2.20 &    2.23 &    0.07 &   -2.79 &   -0.14 \\  -7.98  & -3.28 & -12.94  &  -0.12  &  5.37  &  0.23      \end{pmatrix}$}  &      \textbf{-0.77}          & \multicolumn{1}{l|}{89}     & \multicolumn{1}{l|}{40.5}    & \multicolumn{1}{l|}{  $\begin{pmatrix} -2.12  &  3.87 &   1.47 &   -0.26 &  -0.04 &   0.65 \\  3.75 &  -14.55 &   -1.48 &    1.14 &    5.35 &    2.03\\ -0.67 &    2.20 &    2.23 &    0.07 &   -2.79 &   -0.14 \\  -7.98  & -3.28 & -12.94  &  -0.12  &  5.37  &  0.23      \end{pmatrix}$ }  &     \textbf{-0.77}         \\ \hline
hf2d\_is5  ($n_x = 5$)           & \multicolumn{1}{l|}{\textbf{5}}     & \multicolumn{1}{l|}{\textbf{1.4}}    & \multicolumn{1}{l|}{$\begin{pmatrix} 5.8  &  2.7 &  0.08 &  -0.28 \\  -1.18  & -1.07  &  1.41  &  2.04  \end{pmatrix}$}  &   \textbf{-5.17}           & \multicolumn{1}{l|}{14}     & \multicolumn{1}{l|}{3.5}    & \multicolumn{1}{l|}{$\begin{pmatrix}  5.8  &  2.7 &  0.08 &  -0.28 \\  -1.18  & -1.07  &  1.41  &  2.04 \end{pmatrix}$}  &      \textbf{-5.17}         \\ \hline
hf2d\_cd4  ($n_x = 7$)               & \multicolumn{1}{l|}{\textbf{6}}     & \multicolumn{1}{l|}{\textbf{1.8}}    & \multicolumn{1}{l|}{$\begin{pmatrix} -3.2 &  -3.7 \\ -3.5 & -3.9 \end{pmatrix}$}  &     \textbf{-2.5}            & \multicolumn{1}{l|}{78}     & \multicolumn{1}{l|}{17}    & \multicolumn{1}{l|}{$\begin{pmatrix} -3.3 &  -4.3 \\ -4.3 & -5.5 \end{pmatrix}$}  &      -2.48   \\ \hline
hf2d\_cd5 ($n_x = 7$)                 & \multicolumn{1}{l|}{\textbf{8} }    & \multicolumn{1}{l|}{\textbf{2.5}}    & \multicolumn{1}{l|}{$\begin{pmatrix} -0.23 & -0.21  \\ -1.38 &  -0.43 \end{pmatrix}$}  &  \textbf{-1.79}        & \multicolumn{1}{l|}{257}     & \multicolumn{1}{l|}{19.6}    & \multicolumn{1}{l|}{$\begin{pmatrix} -0.57  & -1.54  \\  -1.54   & -3.65  \end{pmatrix}$}  &       -1.37      \\ \hline
je2     ($n_x = 21$)             & \multicolumn{1}{l|}{\textbf{16}}     & \multicolumn{1}{l|}{\textbf{11.6}}    & \multicolumn{1}{l|}{$\begin{pmatrix}  1.328  &  0.087  & -0.090 \\ -1.462 &   0.1918  &  1.927 \\  1.893  &  0.4696  &  2.7049  \end{pmatrix}$}  &       \textbf{-2.51}           & \multicolumn{1}{l|}{47}     & \multicolumn{1}{l|}{56.5}    & \multicolumn{1}{l|}{$\begin{pmatrix}  1.328  &  0.087  & -0.090 \\ -1.462 &   0.1918  &  1.927 \\  1.893  &  0.4696  &  2.7049  \end{pmatrix}$}  &  \textbf{-2.51}             \\ \hline
lah      ($n_x = 48$)           & \multicolumn{1}{l|}{\textbf{5}}     & \multicolumn{1}{l|}{\textbf{51}}    & \multicolumn{1}{l|}{-6}  &             \textbf{-2.69}    & \multicolumn{1}{l|}{99}     & \multicolumn{1}{l|}{1037.3}    & \multicolumn{1}{l|}{-5.99}  &     \textbf{-2.69}         \\ \hline
rea1     ($n_x = 4$)            & \multicolumn{1}{l|}{\textbf{5}}     & \multicolumn{1}{l|}{\textbf{1.4}}    & \multicolumn{1}{l|}{ $\begin{pmatrix}  -1.740  &  4.229 &  -2.175 \\ 5.147 &  -16.347 &    6.728   \end{pmatrix}$}  &    \textbf{-3}             & \multicolumn{1}{l|}{22}     & \multicolumn{1}{l|}{3.5}    & \multicolumn{1}{l|}{ $\begin{pmatrix}  -1.740  &  4.229 &  -2.175 \\ 5.147 &  -16.347 &    6.728   \end{pmatrix}$}  &      \textbf{-3}         \\ \hline
wec2   ($n_x = 10$)              & \multicolumn{1}{l|}{\textbf{14} }    & \multicolumn{1}{l|}{\textbf{8.5}}    & \multicolumn{1}{l|}{$\begin{pmatrix}-1.0733  &     -0.34109 & 0.9588   &   0.0988   \\ 0.1757 &       -0.1420   &  -1.39116  &       -0.10933   \\ 0.9581    &     0.80115 & 0.19483  &      0.66336   \end{pmatrix}$}  &     \textbf{-1.3796}            & \multicolumn{1}{l|}{40}     & \multicolumn{1}{l|}{28.9}    & \multicolumn{1}{l|}{$\begin{pmatrix} 0.2788  &  0.09640 &  34.9399  &   0.06837 \\  -0.3283    &     -0.1234  & -0.07736  &    -0.00402 \\0.0329    &      8.6410  & 1.3824   &      0.79048 \end{pmatrix}$}  &       -0.6829       \\ \hline
\end{tabular}
\end{adjustbox}
}
\caption{\small Performance of  RHOTA algorithm for $p=1$ and BMIsolver \cite{DiGuMiDi:11} using data from  COMPl$_{e}$ib library \big(MEV = maximum eigenvalue of the real part of $(A + BKC)$\big).}
\label{tab:1}
\vspace{-0.9cm}
\end{table}

\noindent We compare RHOTA algorithm for $p=1$ with  BMIsolver \cite{DiGuMiDi:11}. BMIsolver is specifically designed to optimize the spectral abscissa of the closed-loop system $\dot{x} = (A + BKC)x$ (refer to \cite{DiGuMiDi:11} for comprehensive details). In Table \ref{tab:1}, we report the number of iterations, CPU time, the obtained solution $K$, and the maximum eigenvalue of the real part of the matrix $A + BKC$ (called spectral abscissa). The stopping criterion utilized is $\|F(X_k,K_k,Q_k) \|\leq 10^{-3}$ and we use CVX to solve the subproblem  in RHOTA  \cite{GrBo:13}. Both algorithms, BMIsolver and RHOTA, commence with identical initial values $X_0$ and $K_0$. Each test case from COMPl$_{e}$ib is initialized differently. From the data presented in Table \ref{tab:1}, it is evident that RHOTA  outperforms BMIsolver \cite{DiGuMiDi:11} in terms of both, CPU time and number of iterations. Moreover, RHOTA  yields a slightly smaller value for the spectral abscissa, showcasing the efficiency of the proposed method. The superior performance of RHOTA algorithm (in time and iterations) can be attributed to two facts: first, by linearizing inside the norm, RHOTA leverages a portion of the Hessian of the objective, while BMIsolver solely utilizes first-order information; second, our formulation \eqref{fw: eq3} satisfies the KL property \eqref{eq:kl} (as composition of semi-algebraic functions, i.e., the 2-norm and quadratic functions, see \cite{BolSab:14} and also \cite{FaBo:21}), which, according to Theorem \ref{th:3}, ensures fast convergence for RHOTA.


\subsection{Solving phase retrieval problems}
In this section, we present numerical simulations for  solving the phase retrieval problem \cite{CaLiSo:15,DuRu:19},  using  real images from the collection of handwritten digits, accessible at \cite{MNIS}. The primary objective is to evaluate the performance of  HOPP method in image recovery and compare it with the prox-linear method introduced in \cite{DuRu:19}. Given that \cite{DuRu:19} demonstrates perfect image recovery under real-valued random Gaussian measurements, even when $m=2\times n$, we adopt similar settings. Specifically, we evaluate the performance of our method for $p=1,2$ and the prox-linear method  in \cite{DuRu:19}, aiming to recover a digit image using Gaussian measurement vectors $a_i\in\mathbb{R}^n$ and set $Q_i = a_i^{T}a_i$ for $i=1:m$ with $m = 2\times n$. To initialize the process, we introduce some noise to the real-digit image to generate the starting point $x_0$. The stopping criterion for both methods is set as $f (x_k) \leq 10^{-4}$ or $k\geq 100$. Each subproblem is solved using CVX \cite{GrBo:13}. 

\begin{figure}
    \begin{flushleft}
       \includegraphics[width = 13cm, height = 9cm]{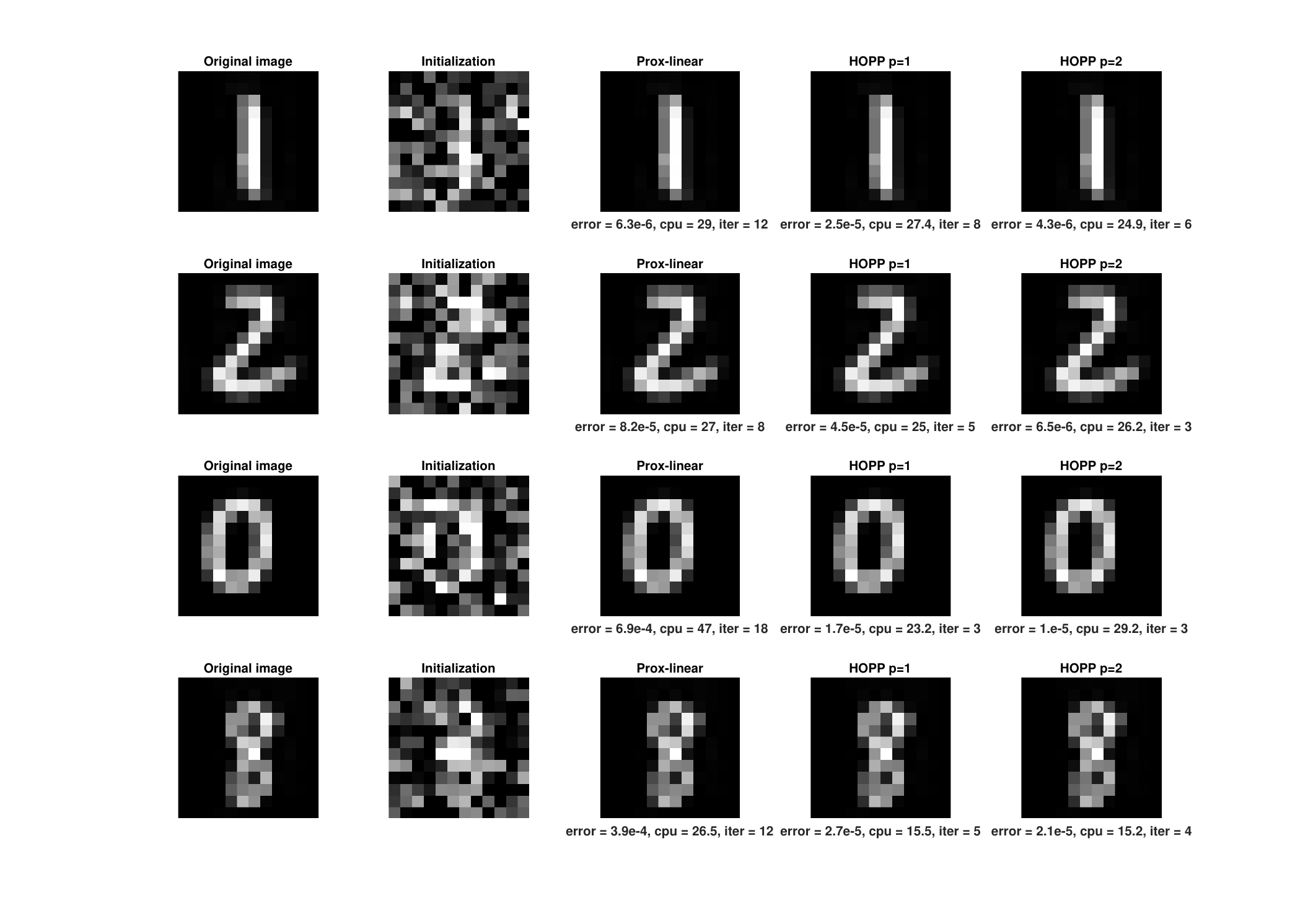}
       \vspace{-1cm}
        \caption{\small Performance of  prox-linear method \cite{DuRu:19} and  HOPP  for $p=1$ and $p=2$ with $M = 0.1$ on  $12\times12$ digit images: initialization satisfying $\text{dist}(x_0,x^*) < \|x^*\| \sigma_0/L$.}
        \label{fig:1}
    \end{flushleft}
    \begin{flushleft}
       \includegraphics[width = 13cm, height = 9cm]{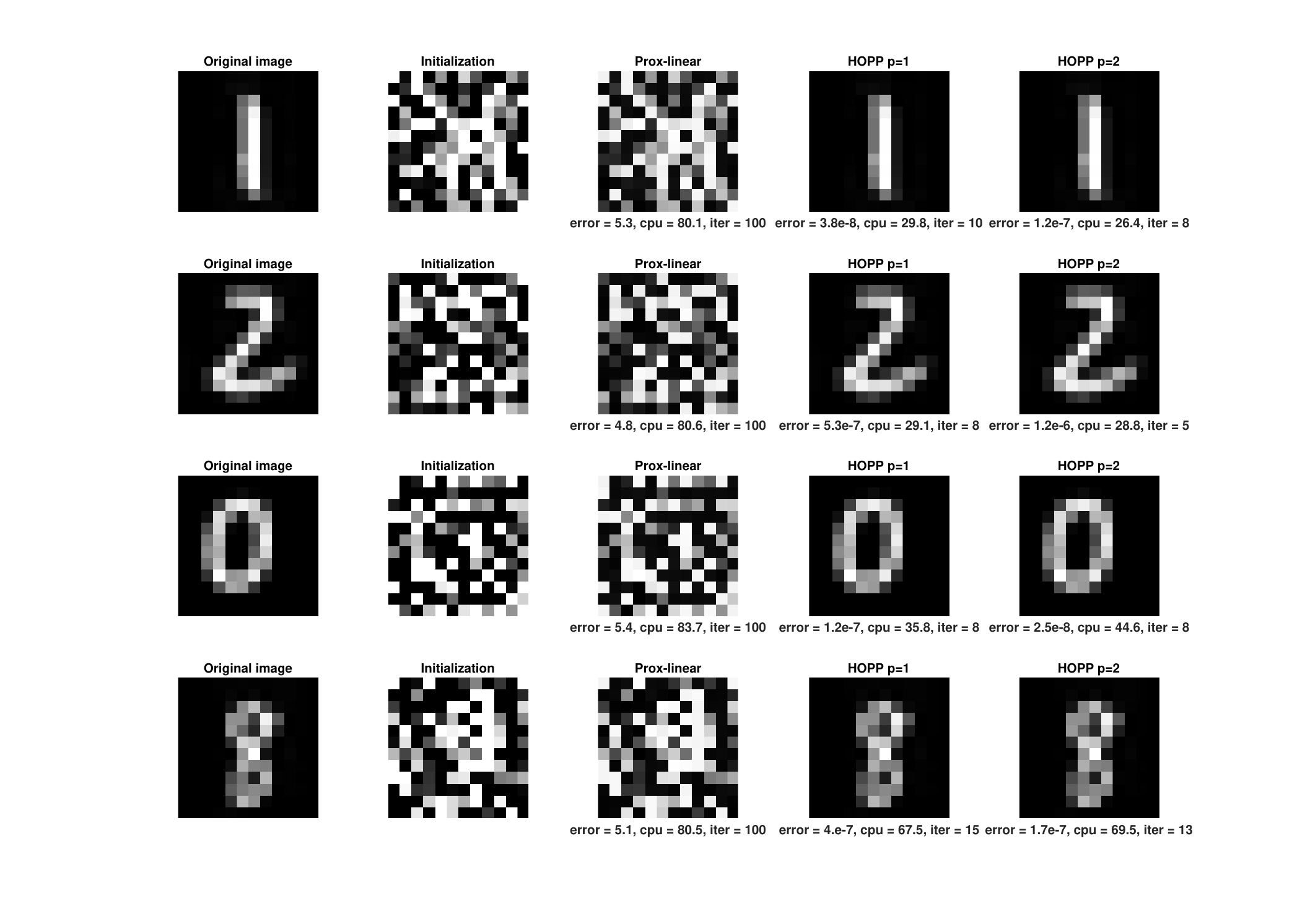}
        \caption{\small Performance of  prox-linear method \cite{DuRu:19} and  HOPP  for $p=1$ and $p=2$, with $M = 0.01$ on $12\times12$ digit images: random initialization $x_0$.}
        \label{fig:2}
    \end{flushleft}
    \vspace{-0.9cm}
\end{figure}

\medskip 

\noindent The results are presented in Figures \ref{fig:1} and \ref{fig:2}. In Figure \ref{fig:1}, we initialize the starting point $x_0$ (by adding some noise to the original image $x^*$) to satisfy the constant relative error guarantee $\|x_0 - x^*\| < \frac{\sigma_0}{L}\|x^*\|$, with $L = \|\left(\|a_1\|^2,\cdots,\|a_m\|^2 \right)\|$, as presented in \cite{DuRu:19}. From Figure \ref{fig:1}, it's evident that both algorithms achieve  good recovery of the original image, with HOPP $(p=2)$ given the best error $\|x_k - x^*\|$. However, HOPP algorithm for $p=1, 2$ is much faster than the prox-linear algorithm \cite{DuRu:19}. In Figure \ref{fig:2}, we set the initial point $x_0$ randomly, so that it does not satisfy the condition $\|x_0 - x^*\| < \frac{\sigma_0}{L}\|x^*\|$. Notably, Figure \ref{fig:2} illustrates that the prox-linear algorithm \cite{DuRu:19} fails to recover the original image after 100 iterations (e.g., the error $\|x_k - x^*\| \approx 5$). In contrast,  HOPP algorithm (for both $p=1$ and $p=2$) is able to recover very well the original image for sufficiently small $M$ (e.g., the error $\|x_k - x^*\| \approx 10^{-6}$). This highlights the efficiency and robustness of  HOPP algorithm. In cases where the true image $x^*$ is unknown, we posit that the HOPP method's flexibility, that follows from  the free choice of the regularization parameter  $M$, allows it to perform effectively even with an initial point that is not necessarily close to the true image. Such an initial point could be generated more affordably than the methods proposed in \cite{DuRu:19,CaLiSo:15}.  Finally,  one can notice from  Figures \ref{fig:1} and \ref{fig:2} the considerable time taken by CVX to solve the convex subproblems. Thus, it would be interesting to explore more efficient convex solvers for solving these subproblems. This aspect remains open for further investigations.






\end{document}